\documentclass[journal,twoside,web]{ieeecolor}
\usepackage{generic}
\usepackage{cite}
\usepackage{amsmath,amssymb,amsfonts}
\usepackage{graphicx}
\usepackage{textcomp}

\usepackage{amsmath,amssymb,amsfonts}
\usepackage{algorithm,algorithmicx,algpseudocode}
\usepackage{color}

\newtheorem{theorem}{\bf Theorem}[section]

\newtheorem{lemma}[theorem]{\bf Lemma}
\newtheorem{remark}[theorem]{\bf Remark}
\newtheorem{assumption}[theorem]{\bf Assumption}

\begin{document}
\title{Structure-Exploiting Newton-Type Method for Optimal Control of Switched Systems}
\author{Sotaro Katayama and Toshiyuki Ohtsuka 
\thanks{Sotaro Katayama and Toshiyuki Ohtsuka are with Graduate School of Informatics, Kyoto University, Kyoto, Japan. 
(e-mail: katayama.25w@st.kyoto-u.ac.jp, ohtsuka@i.kyoto-u.ac.jp).}}

\maketitle

\begin{abstract}
This study proposes an efficient Newton-type method for the optimal control of switched systems under a given mode sequence.
A mesh-refinement-based approach is utilized to discretize continuous-time optimal control problems (OCPs) and formulate a nonlinear program (NLP), which guarantees the local convergence of a Newton-type method.
A dedicated structure-exploiting algorithm (Riccati recursion) is proposed to perform a Newton-type method for the NLP efficiently because its sparsity structure is different from a standard OCP.
The proposed method computes each Newton step with linear time-complexity for the total number of discretization grids as the standard Riccati recursion algorithm. 
Additionally, the computation is always successful if the solution is sufficiently close to a local minimum.
Conversely, general quadratic programming (QP) solvers cannot accomplish this because the Hessian matrix is inherently indefinite.
Moreover, a modification on the reduced Hessian matrix is proposed using the nature of the Riccati recursion algorithm as the dynamic programming for a QP subproblem to enhance the convergence. 
A numerical comparison is conducted with off-the-shelf NLP solvers, which demonstrates that the proposed method is up to two orders of magnitude faster.
Whole-body optimal control of quadrupedal gaits is also demonstrated and shows that the proposed method can achieve the whole-body model predictive control (MPC) of robotic systems with rigid contacts.
\end{abstract}


\section{Introduction}\label{sec:introduction}

\IEEEPARstart{S}{witched} systems are a class of hybrid systems consisting of a finite number of subsystems and switching laws of active subsystems.
Many practical control systems are modeled as switched systems, such as real-world complicated process systems \cite{bib:CIA-MPC}, automotive systems with gear shifts \cite{bib:MIOCPforEV}, and robotic systems with rigid contacts \cite{bib:twoStageApplication1, bib:twoStageApplication2}.
It is difficult to solve the optimal control problems (OCPs) of such systems because these OCPs generally involve mixed-integer nonlinear programs (MINLPs) \cite{bib:MINLP}. 
Therefore, model predictive control (MPC) \cite{bib:mpcbook}, in which OCPs must be solved in real-time, is particularly difficult for switched systems.
One of the most practical approaches for solving the OCPs of switched systems is the combinatorial integral approximation (CIA) decomposition method \cite{bib:CIA, bib:CIA-MPC}. 
CIA decomposition relaxes the MINLP into a nonlinear program (NLP) in which the binary variables are relaxed into continuous variables. Subsequently, the integer variables are approximately reconstructed from the relaxed NLP solution by solving a mixed-integer linear program.
However, vanishing constraints, which typically model mode-dependent path constraints, raise numerical issues due to the violation of the linear independence constraint qualification (LICQ) \cite{bib:VanishingConstraints}.
A similar smoothing approach is used for discrete natures when solving OCPs for mechanical systems with rigid contacts, which are often approximately formulated as mathematical programs with complementarity constraints (MPCCs) \cite{bib:MPCC1, bib:MPCC2}. 
However, this case has the same LICQ problem as the vanishing constraints, requiring a great deal of computational time to alleviate the numerical ill-conditioning. 
Furthermore, it often suffers from undesirable stationary points \cite{bib:MPCCLimits}.

Another tractable and practical approach for the optimal control of switched systems is to fix the switching (mode) sequence.
For example, in robotics applications, a high-level planner computes the feasible contact sequence taking perception into account. Subsequently, the sequence is provided to a lower-level optimal control-based dynamic motion planner or MPC controller \cite{bib:TedrakeAtlas, bib:perceptiveLocomotion, bib:mpc-qp-robot}. 
The lower-level dynamic planner or MPC then discovers the optimal switching times and optimal control input.
An advantage of this approach over the CIA decomposition and MPCC approaches is that the optimization problem is smooth and can therefore be efficiently solved without suffering from the LICQ problem.
For example, \cite{bib:STO1, bib:STO2} proposed efficient Newton-type methods for OCPs of switched systems with autonomous subsystems. That is, the switched systems only included a switching signal but no continuous control input, which is called the switching-time optimization (STO) problem.
The STO approach can be more efficient than the CIA decomposition because the STO problem can be formulated as a smooth and tractable NLP,  as numerically shown in \cite{bib:STO2}.
However, once the switched system includes a continuous control input, the efficient methods of \cite{bib:STO1, bib:STO2} cannot be applied.
To the best of our knowledge, there is no efficient numerical method for OCPs of such systems. As a result, many real-world robotic applications are limited to focusing on dynamic motion planning with fixed switching instants \cite{bib:MPCforQuadruped1, bib:MITCheetahMPC, bib:PlannningTerrainMapping}.

The two-stage approach has been studied for the OCPs of switched systems with continuous control input \cite{bib:two-stage0, bib:twoStage, bib:twoStageApplication1, bib:twoStageApplication2, bib:twoStageConstrained}.
A general two-stage approach was proposed in \cite{bib:two-stage0, bib:twoStage}.
In this approach, the optimization problem to determine the control input and switching instants was decomposed into an STO problem with a fixed control input (upper-level problem) and standard OCP that only determined the control input with fixed switching instants (lower-level problem).
In this framework, off-the-shelf OCP solvers can be used for a lower-level problem. \cite{bib:two-stage0, bib:twoStage} used an indirect OCP solver to compute an accurate solution of a lower level problem.
\cite{bib:twoStageApplication1, bib:twoStageApplication2} formulated a lower-level problem as a direct OCP and solved it using off-the-shelf Newton-type methods.
However, these studies still lack convergence speed because each of these two stages does not take the other stage into account when solving its own optimization problem. 
As a result, the application examples of \cite{bib:twoStageApplication1, bib:twoStageApplication2} were limited to off-line computation for the trajectory optimization problems of simplified robot models.
Moreover, they could not guarantee the local convergence or address inequality constraints.
\cite{bib:twoStageConstrained} was the first study that guaranteed convergence with a finite number of iterations and treated inequality constraints within a two-stage framework.
However, it still required extensive computational time until convergence, even for a very simple linear quadratic example.

Other approaches simultaneously optimized the switching instants and other variables, such as the state and control input \cite{bib:Betts, bib:GPOPS, bib:hybridMPC, bib:hybridOCPRiccati}.
A multi-phase trajectory optimization \cite{bib:Betts, bib:GPOPS} naturally incorporated the STO problem into the direct transcription, solving the NLP to simultaneously determine all variables (including the state, control input, and switching instants) using general-purpose off-the-shelf NLP solvers, such as Ipopt \cite{bib:ipopt}.
However, the computational speed of general-purpose linear solvers used in off-the-shelf NLPs can typically be further improved, particularly for large-scale systems, because they have a certain sparsity structure \cite{bib:RiccatiMPC, bib:MPCBoyd, bib:forces}.
In \cite{bib:hybridMPC}, we applied the simultaneous approach with the direct single-shooting method and achieved real-time MPC for a simple walking robot using the Newton-Krylov type method \cite{bib:cgmres}. 
In \cite{bib:hybridOCPRiccati}, we formulated an NLP using the direct multiple-shooting method \cite{bib:DMS} and proposed a Riccati recursion algorithm for the NLP, which achieved faster computational time compared with the two-stage methods.
However, these methods lacked the convergence guarantee due to the irregular discretization of the continuous-time OCP, which changed the problem structure throughout the Newton-type iterations. 
As a result, the method from \cite{bib:hybridOCPRiccati} could only converge when the initial guess of the switching instants was close to the optimal one in the numerical example.

This study proposes an efficient Newton-type method for optimal control of switched systems under a given mode sequence.
First, a mesh-refinement-based approach is proposed to discretize the continuous-time OCP using the direct multiple-shooting method \cite{bib:DMS} to formulate an NLP that facilitates the local convergence of the Newton-type methods.
Second, a dedicated efficient structure-exploiting algorithm (Riccati recursion algorithm \cite{bib:RiccatiMPC}) is proposed for the Newton-type method because the sparsity structure of the NLP is different from the standard OCP.
The proposed method computes each Newton step with linear time-complexity of the total number of the discretization grids as the standard Riccati recursion algorithm.
Additionally, it can always solve the Karush-Kuhn-Tucker (KKT) systems arising in the Newton-type method if the solution is sufficiently close to a local minimum, so that the second-order sufficient condition (SOSC) holds.
This is in contrast to some general quadratic programming (QP) solvers that cannot treat the proposed formulation because the Hessian matrix is inherently indefinite.
Third, a modification on the reduced Hessian matrix is proposed to enhance the convergence using the nature of the Riccati recursion algorithm as the dynamic programming \cite{bib:DPandOCP} for a QP subproblem. 
Two numerical experiments are conducted to demonstrate the efficiency of the proposed method: a comparison with off-the-shelf NLP solvers and testing the whole-body optimal control of quadrupedal gaits.
The comparison with off-the-shelf NLP solvers showed that the proposed method could solve the OCPs that the sequential quadratic programming (SQP) method with qpOASES \cite{bib:qpOASES} or OSQP \cite{bib:OSQP} failed to solve and was up to two orders of magnitude faster than a general NLP solver (Ipopt) \cite{bib:ipopt}.
The whole-body optimal control of quadrupedal gaits showed that the proposed method achieved the whole-body MPC of robotic systems with rigid contacts.

The remainder of this paper is organized as follows.
A mesh-refinement-based discretization method of the continuous-time OCP is presented in Section \ref{sec:NLP}.
The KKT system to be solved to compute the Newton-step is discussed in Section \ref{sec:NewtonTypeMethod} and the Riccati recursion algorithm to solve the KKT system and its convergence properties are described in Section \ref{sec:Riccati}.
The reduced Hessian modification using the Riccati recursion algorithm is also described in this section.
The above formulations and Riccati recursion algorithm are extended to switched systems with state jumps and switching conditions in Section \ref{sec:Riccati}, representing robotic systems with contacts.
A numerical comparison of the proposed method with off-the-shelf NLP solvers and examples of whole-body optimal control of a quadrupedal robot are presented in Section \ref{sec:experiments}. 
Finally, a brief summary and mention of future work are presented in Section \ref{sec:conclu}.

\noindent
\textit{Notation and preliminaries: } The Jacobians and Hessians of a differentiable function using certain vectors are described as follows: $\nabla_x f (x)$ denotes $\left( \frac{\partial f}{\partial x} \right) ^{\rm T} (x)$, and $\nabla_{x y} g(x, y)$ denotes $\frac{\partial^2 g}{\partial x \partial y} (x, y)$.
A diagonal matrix whose elements are a vector $x$ is denoted as ${\rm diag} (x)$.
A vector of an appropriate size with all components presented by $\alpha \in \mathbb{R}$ is denoted as $\alpha 1$.
All functions are assumed to be twice-differentiable.

\section{Problem Formulation}\label{sec:NLP}

\begin{figure}[tb]
\centering
\includegraphics[scale=0.40]{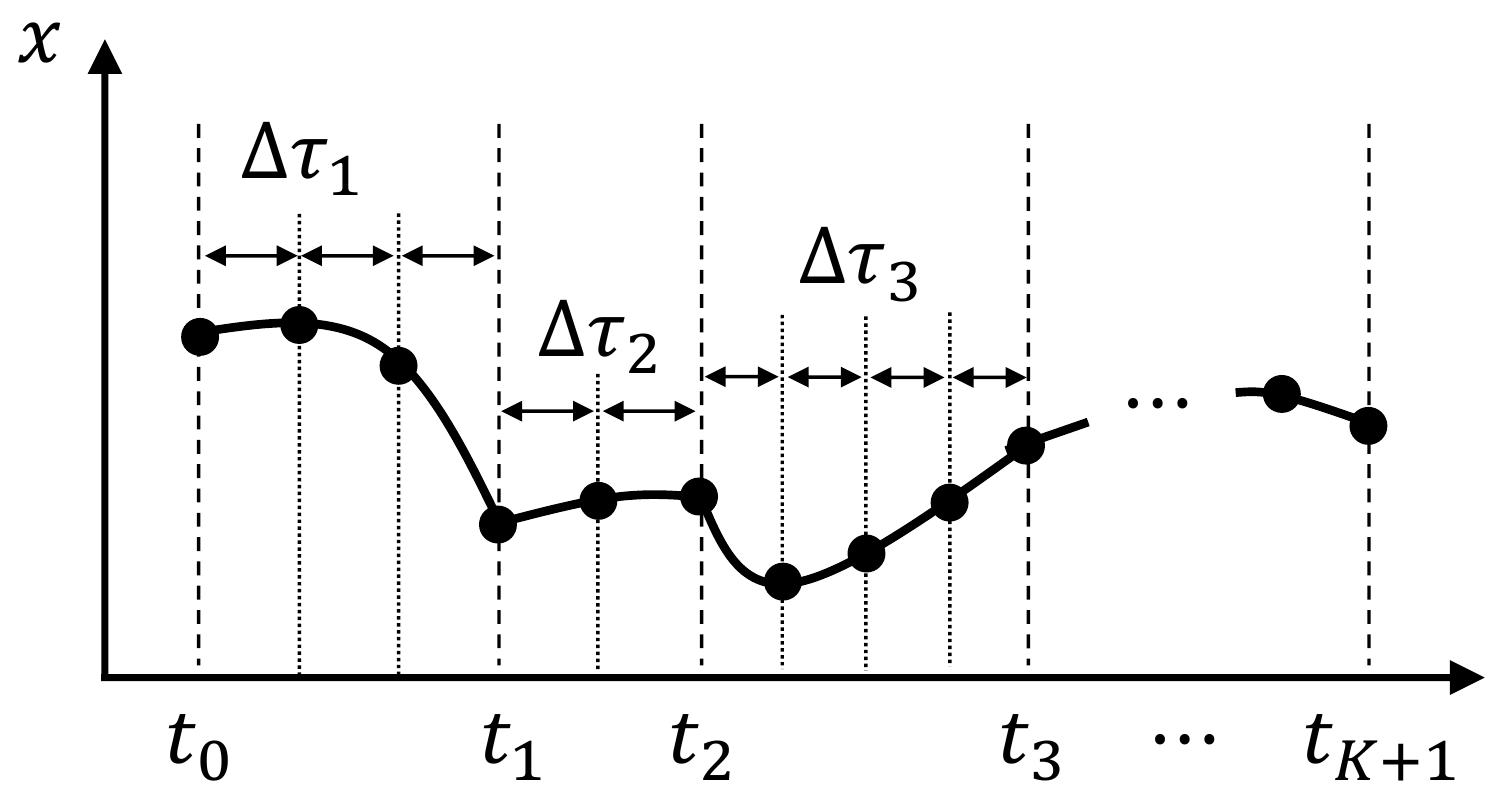}
\caption{Proposed discretization method for the OCPs of switched systems.}
\label{fig:discretization}
\end{figure}

We consider a switched system consisting of $K+1$ ($ K > 0$) subsystems, which is expressed as: 
\begin{equation}\label{eq:stateEquation}
    \dot{x} (t) = f_{k} (x(t), u(t)), \;\; t \in [t_{k-1}, t_{k}], \;\; k \in \mathcal{K}.
\end{equation}
The system also contains constraints, which are expressed as: 
\begin{equation}\label{eq:constraints}
    g_{k} (x(t), u(t)) \leq 0, \;\; t \in [t_{k-1}, t_{k}], \;\; k \in \mathcal{K} ,
\end{equation}
where $x (t) \in \mathbb{R}^{n_x}$ denotes the state, $u (t) \in \mathbb{R}^{n_u}$ denotes the control input, $f_k : \mathbb{R}^{n_x} \times \mathbb{R}^{n_u} \to \mathbb{R}^{n_x}$, and $g_k : \mathbb{R}^{n_x} \times \mathbb{R}^{n_u} \to \mathbb{R}^{n_g}$. $\mathcal{K} := \left\{ 1, ..., K+1 \right\}$ denotes the given indices of the active subsystems. 
Note that the index $k$ is also referred to as a ``phase" in this study.
$t_0$ and $t_{K+1}$ denote the fixed initial and terminal times of the horizon, respectively. $t_k$, $k \in \left\{ 1, ..., K \right\}$ denotes the switching instant from phase $k$ to phase $k+1$.
The OCP of the switched system for a given initial state $x(t_0) \in \mathbb{R}^{n_x}$ is expressed as: 
\begin{subequations}\label{subeq:continuousTimeOCP}
\begin{equation}\label{eq:costFunction}
    \min_{u(\cdot), t_{1}, ..., t_{K}} {J} = \; V_f (x(t_{K+1})) + \sum_{k=1}^{K+1} \int_{t_{k-1}}^{t_{k}} l_k (x(\tau), u(\tau)) d \tau
\end{equation}
\begin{align}\label{eq:tkConstraints}
    {\rm s.t.} \;\;\; & {\rm (\ref{eq:stateEquation})}, \; {\rm (\ref{eq:constraints})}, \; \notag \\ 
    & t_{k-1} + \underline{\Delta}_{k} \leq t_{k}, \;\; k \in \left\{ 1, ..., K \right\},
\end{align}
\end{subequations}
where $\underline{\Delta}_{k} \geq 0, \; k \in \mathcal{K}$ is the minimum dwell-time.

Next, a discretization method and mesh-refinement-based solution approach are proposed for the OCP to guarantee the local convergence.
The OCP is discretized with the direct multiple-shooting method \cite{bib:DMS} based on the forward Euler method so that all the time-steps in phase $k \in \mathcal{K}$ are equal.
Note that it is easy to extend the proposed method for higher-order explicit integration schemes to the direct multiple-shooting method, such as the fourth-order explicit Runge-Kutta method, as long as the Riccati recursion algorithm can be applied to the integration schemes.
We introduce $N$ grid points over the horizon, the discretized state $X := \left\{x_0, ..., x_N \right\}$, discretized control input $U := \left\{ u_0, ..., u_{N-1} \right\}$, and switching instants $T:= \left\{t_1, ..., t_K \right\}$.
The NLP is then expressed as:
\begin{subequations}\label{subeq:NLP}
\begin{equation}\label{eq:NLP:cost}
    \min_{X, U, T} J = V_f (x_N) + \sum_{k \in \mathcal{K}} \sum_{i \in \mathcal{I}_k} l_k (x_i, u_i) \Delta \tau_k
\end{equation}
\begin{align}
    {\rm s.t.} \;\;\; & x_0 - \bar{x} = 0 , \label{eq:NLP:initialStateConstraint} \\
    & x_i + f_k (x_i, u_i) \Delta \tau_k - x_{i+1} = 0,  \;\; i \in \mathcal{I}_k, \; k \in \mathcal{K} , \label{eq:NLP:stateEquation} \\ 
    & g_k (x_i, u_i) \leq 0,  \;\; i \in \mathcal{I}_k, \; k \in \mathcal{K}, \label{eq:NLP:inequalityConstraints} \\
    & t_{k-1} + \underline{\Delta}_k - t_{k} \leq 0, \;\; k \in \left\{ 1, ..., K \right\} \label{eq:NLP:STOConstraints},
\end{align}
where $\mathcal{I}_k$ is the set of stage indices at phase $k$ (that is, the set of time stages where the subsystem equation $f_k (x_i, u_i)$ is active). $\Delta \tau_k$ is the time step at phase $k$, defined as:
\begin{equation}
    \Delta \tau_k := \frac{t_{k} - t_{k-1}}{N_k}, \;\; k \in \left\{ 1, ..., K \right\}.
\end{equation}
\end{subequations}
Note that the last stage $N$ is not included in any $\mathcal{I}_k$ and that the initial state condition is lifted as (\ref{eq:NLP:initialStateConstraint}) for an efficient MPC implementation \cite{bib:RTI}.

The continuous-time OCP (\ref{subeq:continuousTimeOCP}) is solved using an adaptive mesh-refinement approach \cite{bib:Betts}, which consists of a solution of the NLP (\ref{subeq:NLP}) using a Newton-type method and mesh-refinement due to the changes of $\Delta \tau_k$, as shown in Algorithm \ref{alg1}.
After solving the NLP, the size of the discretization steps $\Delta \tau_k$ is checked for each $k \in \mathcal{K}$.
If the step is too large, that is, if it exceeds the specified threshold $\Delta \tau_{\rm max}$, the solution at phase $k$ may be not accurate. 
Therefore, the number of grids for phase $k$ is increased.
Conversely, if the step is too small ($\Delta \tau < \Delta \tau_{\rm min}$), the number of grids for phase $k$ is reduced to decrease the computational time of the next NLP step.
The algorithm terminates if some criteria (for example, $l_2$ norm of the KKT residual), which is denoted as "NLP error" in Algorithm \ref{alg1}, is smaller than a predefined threshold $\epsilon > 0$ and the discretization steps pass the checks.
By appropriately choosing $\Delta \tau_{\rm min}$, the NLP dimension is almost constant throughout Algorithm \ref{alg1}. This is an advantage of the proposed method over direct transcription methods, such as \cite{bib:GPOPS}, whose NLP dimension is unknown before being solved.

\begin{algorithm}[tb]
\caption{Adaptive mesh-refinement approach for the optimal control of switched systems (\ref{subeq:continuousTimeOCP})}
\label{alg1}
\begin{algorithmic}[1]
    \Require The initial state ${x} (t_0)$, initial guess of the solution $X$, $U$, $T$, initial guess of the Lagrange multipliers, maximum and minimum discretization step sizes $\Delta \tau_{\rm max} > \Delta \tau_{\rm min} > 0$, and convergence tolerance $\epsilon > 0$.
    \Ensure Optimal solution $X$, $U$, $T$.
    \While{NLP error $> \epsilon$}
        \If{$\Delta \tau_k < \Delta \tau_{\rm min}$} 
            \State Mesh-refinement for phase $k$ (decrease grids).
        \EndIf
        \State Solve the NLP (\ref{subeq:NLP}) using the Newton-type method.
        \If{$\Delta \tau_k > \Delta \tau_{\rm max}$} 
            \State Mesh-refinement for phase $k$ (increase grids).
        \EndIf
    \EndWhile
\end{algorithmic}
\end{algorithm}

It is worth noting that the NLP (\ref{subeq:NLP}) is smooth and its structure does not change within each NLP step.
Therefore, the Newton-type method for the NLP (\ref{subeq:NLP}) is always tractable.
\begin{remark}\label{remark:localConvergence}
It is trivial to show the local convergence of Newton-type methods for the NLP (\ref{subeq:NLP}) under some reasonable assumptions \cite{bib:nocedal}.
Moreover, if the solution guess after the mesh-refinement of Algorithm \ref{alg1} is sufficiently close to a local minimum of the NLP after the mesh-refinement (that is, the mesh-refinement is sufficiently accurate), then the local convergence of the overall Algorithm \ref{alg1} is also guaranteed.
In contrast, previous methods cannot guarantee the local convergence, such as the two-stage approaches \cite{bib:two-stage0, bib:twoStage, bib:twoStageApplication1, bib:twoStageApplication2} and simultaneous approaches with the other discretization methods \cite{bib:hybridMPC, bib:hybridOCPRiccati}.
\end{remark}

Note that it is assumed throughout this study that the state equations (\ref{eq:stateEquation}), inequality constraints (\ref{eq:constraints}), and cost function (\ref{eq:costFunction}) do not depend on the time. That is, they are time invariant for notational simplicity.
However, it is easy to extend the proposed method to time-varying cases, where the time of each grid $i$ that depends on the switching times $t_1, ..., t_K$ in (\ref{eq:stateEquation}), (\ref{eq:constraints}), and (\ref{eq:costFunction}) must be taken into account.
Therefore, additional sensitivities of (\ref{eq:stateEquation}), (\ref{eq:constraints}), and (\ref{eq:costFunction}) with respect to the switching times are introduced in the KKT conditions and systems, which are derived in the next section.  
The formulations and methods of this study can be directly applied because such additional sensitivities do not change the NLP structure.

\section{KKT System for Newton-Type Method}\label{sec:NewtonTypeMethod}
Next, the KKT system is derived to compute the Newton step of the NLP (\ref{subeq:NLP}).
The inequality constraints are treated with the primal-dual interior point method \cite{bib:nocedal, bib:ipopt}. That is, the slack variables $z_0, ..., z_{N-1} \in \mathbb{R}^{n_g}$ and $w_1, ..., w_K \in \mathbb{R}$ are introduced for (\ref{eq:NLP:inequalityConstraints}) and (\ref{eq:NLP:STOConstraints}), respectively. 
The equality constraints are then considered, which are expressed as:
\begin{subequations}
\begin{align}
    & r_{g, i} := g_k (x_i, u_i) + z_i = 0,  \;\; i \in \mathcal{I}_k, \; k \in \mathcal{K}, \label{eq:PDIPM:inequalityConstraints} \\
    & r_{\Delta, k} := t_{k-1} + \underline{\Delta}_k - t_{k} + w_k = 0, \;\; k \in \left\{ 1, ..., K \right\} 
    \label{eq:PDIPM:STOConstraints} 
\end{align}
\end{subequations}
instead of the inequality constraints (\ref{eq:NLP:inequalityConstraints}) and (\ref{eq:NLP:STOConstraints}).
The barrier functions $ - \epsilon \sum_{i=0}^{N-1} \ln{z_i} - \epsilon \sum_{k=1}^{K} \ln{w_k} $ are also added to the cost function (\ref{eq:NLP:cost}), where $\epsilon > 0$ is the barrier parameter.
The perturbed KKT conditions \cite{bib:nocedal} are obtained by introducing the Lagrange multipliers. 
$\lambda_0, ..., \lambda_N \in \mathbb{R}^{n_x}$ are introduced as the Lagrange multipliers with respect to (\ref{eq:NLP:initialStateConstraint}) and (\ref{eq:NLP:stateEquation}), $\nu_0, ..., \nu_{N-1} \in \mathbb{R}^{n_g}$ are with respect to (\ref{eq:PDIPM:inequalityConstraints}), and $\upsilon_0, ..., \upsilon_{K+1} \in \mathbb{R}$ are with respect to (\ref{eq:PDIPM:STOConstraints}). 
The perturbed KKT conditions are expressed as:
\begin{subequations}\label{subeq:KKTconditions}
\begin{equation}\label{eq:KKTTerminal}
    r_{x, N} := \nabla_x V_f (x_N) - \lambda_N = 0,
\end{equation}
\begin{align}\label{eq:KKTHx}
    r_{x, i} := \; & \nabla_{x} H_k (x_i, u_i, \lambda_{i+1}) \Delta \tau_k + \nabla_x g_k ^{\rm T} (x_i, u_i) \nu_i \notag \\ 
    & + \lambda_{i+1} - \lambda_{i} = 0, \;\;\;
    i \in \mathcal{I}_k, \, k \in \mathcal{K} , 
\end{align}
\begin{align}\label{eq:KKTHu}
    {r}_{u, i} := \nabla_{u} H_k (x_i, u_i, \lambda_{i+1}) \Delta \tau_k + & \nabla_u g_k ^{\rm T} (x_i, u_i) \nu_i = 0, \notag \\  
    & i \in \mathcal{I}_k, \, k \in \mathcal{K} ,
\end{align}
\begin{equation}\label{eq:KKTpdipmPathConstraint} 
    r_{z, i} := {\rm diag} (z_i) \nu_i - \epsilon 1 = 0,  \;\;\; i \in \mathcal{I}_k,  \; k \in \mathcal{K}, 
\end{equation}
\begin{equation}\label{eq:KKTpdipmSTOConstraint}
    r_{w, k} := w_k \upsilon_k - \epsilon = 0, \;\;\; k \in \mathcal{K},
\end{equation}
and
\begin{align}\label{eq:KKTSTO}
    & \frac{1}{N_k} \sum_{i \in \mathcal{I}_k} H_k (x_i, u_i, \lambda_{i+1}) - \frac{1}{N_{k+1}} \sum_{i \in \mathcal{I}_{k+1}} H_{k+1} (x_i, u_i, \lambda_{i+1}) \notag \\ 
    & + \upsilon_k - \upsilon_{k+1} = 0 , \;\;\; k \in \left\{ 2, ..., K+1 \right\} ,
\end{align}
where 
\begin{equation}
    H_k (x_i, u_i, \lambda_{i+1}) := l_k (x_i, u_i) + \lambda_{i+1} ^{\rm T} f_k (x_i, u_i) 
\end{equation}
is the Hamiltonian at phase $k \in \mathcal{K}$.
\end{subequations}
Next, the KKT system is derived to compute the Newton steps of all variables, that is,
$\Delta x_0, ..., \Delta x_N \in \mathbb{R}^{n_x}$, $\Delta u_0, ..., \Delta u_{N-1} \in \mathbb{R}^{n_u}$ $\Delta t_1, ..., \Delta t_{K} \in \mathbb{R}$, and $\Delta \lambda_0, ..., \Delta \lambda_N \in \mathbb{R}^{n_x}$.
Note that the KKT system is herein considered as the standard primal-dual interior point method \cite{bib:nocedal, bib:ipopt}, in which the Newton directions regarding the inequality constraints are eliminated (that is, $\Delta z_0, ..., z_{N-1}, \Delta \nu_0, ..., \Delta \nu_{N-1} \in \mathbb{R}^{n_g}$ and $\Delta w_0, ..., \Delta w_{K+1}, \Delta \upsilon_0, ..., \Delta \upsilon_{K+1} \in \mathbb{R}$). 
The KKT system of interest is then expressed as: 
\begin{subequations}\label{subeq:KKTsystem}
\begin{equation}\label{eq:initialStateLinearized}
    \Delta x_0 + x_0 - \bar{x} = 0,
\end{equation}
\begin{align}\label{eq:stateEquationLinearized}
    A_i \Delta x_i + B_i \Delta u_i + f_i (\Delta t_k - \Delta t_{k-1}) & - \Delta x_{i+1} + \bar{x}_i = 0, \notag \\
    & i \in \mathcal{I}_k, \;\; k \in \mathcal{K},
\end{align}
\begin{align}\label{eq:LxLinearized}
    & Q_{xx, i} \Delta x_i + Q_{xu, i} \Delta u_i + A_i ^{\rm T} \Delta \lambda_{i+1} - \Delta \lambda_{i} \notag \\ 
    & + h_{x, i} (\Delta t_k - \Delta t_{k-1}) + \bar{l}_{x, i} = 0, \;\; i \in \mathcal{I}_k, \;\; k \in \mathcal{K},
\end{align}
\begin{align}\label{eq:LuLinearized}
    & Q_{xu, i} ^{\rm T} \Delta x_i + Q_{uu, i} \Delta u_i + B_i ^{\rm T} \Delta \lambda_{i+1} \notag \\
    & + h_{u, i} (\Delta t_k - \Delta t_{k-1}) + \bar{l}_{u, i} = 0, \;\; i \in \mathcal{I}_k, \;\; k \in \mathcal{K},
\end{align}
\begin{align}\label{eq:STOLinearized}
    & \frac{1}{N_k} \sum_{i \in \mathcal{I}_k} \left( h_{x, i} ^{\rm T} \Delta x_i + h_{u, i} ^{\rm T} \Delta u_i + f_i ^{\rm T} \Delta \lambda_{i+1} + \bar{h}_i \right) \notag \\ 
    & - \frac{1}{N_{k+1}} \sum_{i \in \mathcal{I}_{k+1}} \left( h_{x, i} ^{\rm T} \Delta x_i + h_{u, i} ^{\rm T} \Delta u_i + f_i ^{\rm T} \Delta \lambda_{i+1} + \bar{h}_i \right) \notag \\ 
    & + Q_{tt, k} (\Delta t_{k} - \Delta t_{k-1}) + \bar{q}_{t, k} = 0, \;\;  k \in \left\{ 1, ..., K \right\} ,
\end{align}
and
\begin{equation}\label{eq:KKTSystemTerminal}
    Q_{xx, N} \Delta x_N - \Delta \lambda_{N} + \bar{l}_{x, N} = 0,
\end{equation}
where 
\begin{align*}
    Q_{xx, i} := \; & \nabla_{x x} H_k (x_i, u_i, \lambda_{i+1}) \Delta \tau_k \notag \\ 
    & + \nabla_{x} g_k ^{\rm T} (x_i, u_i) {\rm diag}(z_i) ^{-1} {\rm diag}(\nu_i) \nabla_{x} g_k ^{\rm T} (x_i, u_i),
\end{align*}
\begin{align*}
    Q_{xu, i} := \; & \nabla_{x u} H_k (x_i, u_i, \lambda_{i+1}) \Delta \tau_k \notag \\ 
    & + \nabla_{x} g_k ^{\rm T} (x_i, u_i) {\rm diag}(z_i) ^{-1} {\rm diag}(\nu_i) \nabla_{u} g_k ^{\rm T} (x_i, u_i), 
\end{align*}
\begin{align*}
    Q_{uu, i} := \; & \nabla_{u u} H_k (x_i, u_i, \lambda_{i+1}) \Delta \tau_k \notag \\ 
    & + \nabla_{u} g_k ^{\rm T} (x_i, u_i) {\rm diag}(z_i) ^{-1} {\rm diag}(\nu_i) \nabla_{u} g_k ^{\rm T} (x_i, u_i), 
\end{align*}
\begin{equation*}
    \bar{l}_{x, i} := r_{x, i} + \nabla_{x} g_k ^{\rm T} (x_i, u_i) {\rm diag}(z_i) ^{-1} ({\rm diag}(\nu_i) r_{g, i} - r_{z, i}) ,
\end{equation*}
\begin{equation*}
    \bar{l}_{u, i} := r_{u, i} + \nabla_{u} g_k ^{\rm T} (x_i, u_i) {\rm diag}(z_i) ^{-1} ({\rm diag}(\nu_i) r_{g, i} - r_{z, i}) , 
\end{equation*}
\begin{equation*}
    Q_{tt, k} := w_k ^{-1} \upsilon_k, \;\;
    \bar{q}_{t, k} := w_k ^{-1} (\upsilon_k r_{\Delta, k} - r_{w, k}),
\end{equation*}
\begin{equation*}
    A_i := I + \nabla_x f_k (x_i, u_i) \Delta \tau_k, \;\;
    B_i := \nabla_x f_k (x_i, u_i) \Delta \tau_k, 
\end{equation*}
\begin{equation*}
    f_i := \frac{1}{N_k} f_k (x_i, u_i), \;\;\;  h_{i} := \frac{1}{N_k} H_k (x_i, u_i, \lambda_i), 
\end{equation*}
and
\begin{equation*}
    h_{x, i} := \frac{1}{N_k} \nabla_x H_k (x_i, u_i, \lambda_i), \;\;\; 
    h_{u, i} := \frac{1}{N_k} \nabla_u H_k (x_i, u_i, \lambda_i).
\end{equation*}
\end{subequations}
In addition, $\bar{x}_i$ is residual in (\ref{eq:NLP:stateEquation}). 
Note that some of the phase index $k$ is omitted from the KKT system equations (\ref{subeq:KKTsystem}) because the matrices and vectors (other than the Newton steps in (\ref{subeq:KKTsystem})) are fixed once they are computed and therefore do not depend on the phase index $k$ to solve the KKT system.

\begin{figure*}
\begin{align}\label{eq:QPsubproblem}
     & \min_{\substack{\Delta u_0, ..., \Delta u_{N-1} \\ \Delta t_1, ..., \Delta t_{K+1}}} \sum_{k \in \mathcal{K}} \; \sum_{i \in \mathcal{I}_k,  i \neq N} \frac{1}{2} \left\{
    \begin{bmatrix}
        \Delta t_k - \Delta t_{k-1}  \\ \Delta x_i \\ \Delta u_i 
    \end{bmatrix}^{\rm T}
    \begin{bmatrix}
    0 & h_{x, i} ^{\rm T} & h_{u, i} ^{\rm T} \\
    h_{x, i}  & Q_{xx, i} & Q_{xu, i} \\
    h_{u, i}  & Q_{ux, i} & Q_{uu, i} \\
    \end{bmatrix}
    \begin{bmatrix}
        \Delta t_k - \Delta t_{k-1}  \\ \Delta x_i \\ \Delta u_i 
    \end{bmatrix}
    + \begin{bmatrix}
        \bar{h}_i \\ \bar{l}_{x, i} \\ \bar{l}_{u, i}
    \end{bmatrix}^{\rm T}
    \begin{bmatrix}
        \Delta t_k - \Delta t_{k-1}  \\ \Delta x_i \\ \Delta u_i 
    \end{bmatrix} \right\} \notag \\ 
    & \;\;\;\;\;\;\;\;\;\;\;\;\;\;\;\;\;\;\;\; 
    + \sum_{k \in \mathcal{K}} \left\{ \frac{1}{2} Q_{tt, k} \Delta t_k ^2 - Q_{tt, k} \Delta t_k \Delta t_{k-1} + \bar{q}_{t, k} \Delta t_k \right\} + \Delta x_N ^{\rm T} Q_{xx, N} \Delta x_N + \bar{l}_{x, N} ^{\rm T} \Delta x_N \notag \\
    & \;\;\;\;\;\;\;\; {\rm s.t.}  \; (\ref{eq:initialStateLinearized}), (\ref{eq:stateEquationLinearized}).
\end{align}
\end{figure*}

Note that the KKT system (\ref{subeq:KKTsystem}) is equivalent to the KKT conditions of a QP subproblem (\ref{eq:QPsubproblem}), which is a quadratic approximation of the NLP (\ref{subeq:NLP}).
Subsequently, $\Delta \lambda_0, ..., \Delta \lambda_N$ can be regarded as the Lagrange multiplier of the QP with respect to (\ref{eq:initialStateLinearized}) and (\ref{eq:stateEquationLinearized}) \cite{bib:nocedal}.

\begin{remark}\label{remark:indefiniteHessian}
The Hessian matrix of (\ref{eq:QPsubproblem}) is inherently indefinite, which makes solving the KKT system (\ref{subeq:KKTsystem}) difficult when using off-the-shelf QP solvers because they typically require a positive definite Hessian matrix.
This can be explained with a block diagonal of the Hessian matrix (\ref{eq:QPsubproblem}), expressed as: 
\begin{equation}
    \begin{bmatrix}
        0        & h_{x, i} ^{\rm T} & h_{u, i} ^{\rm T} \\
        h_{x, i} & Q_{xx, i} & Q_{uu, i} \\ 
        h_{u, i} & Q_{xu, i} ^{\rm T} & Q_{uu, i} 
    \end{bmatrix}.
\end{equation}
This is indefinite due to the off-diagonal terms $h_{x, i}$ and $h_{u, i}$, even when: 
$\begin{bmatrix}
    Q_{xx, i} & Q_{uu, i} \\ 
    Q_{xu, i} ^{\rm T} & Q_{uu, i} 
\end{bmatrix} \succ O$.
\end{remark}

\section{Riccati Recursion to Solve KKT Systems}\label{sec:Riccati}
In this section, a Riccati recursion algorithm is presented to compute the Newton step of the NLP (\ref{subeq:NLP}) by solving the KKT system (\ref{subeq:KKTsystem}).
The sparsity structure of the KKT system (\ref{subeq:KKTsystem}) is no longer the same as the standard OCP, which prevents applying the off-the-shelf efficient Newton-type algorithms for OCPs  \cite{bib:RiccatiMPC, bib:MPCBoyd, bib:forces}.
Moreover, as stated in Remark \ref{remark:indefiniteHessian}, the Hessian matrix of the KKT system is indefinite and unsolvable using general QP solvers.
Motivated by these problems, we propose a Riccati recursion algorithm that efficiently solves the KKT system, specifically with $O(N)$ computational time, whenever the SOSC holds.
Moreover, a reduced Hessian modification method is proposed to enhance the convergence when the SOSC does not hold, that is, when the reduced Hessian matrix is indefinite.
As the standard Riccati recursion \cite{bib:RiccatiMPC}, the proposed method is composed of backward and forward recursions, which recursively eliminate the variables from the KKT system (\ref{subeq:KKTsystem}) backward in time and recursively compute the Newton step forward in time, respectively.

\subsection{Backward Recursion}\label{subsec:backwardRiccatiRecursion}
In the backward recursion, the Newton steps are recursively eliminated from stages $N$ to $0$.
Specifically, expressions of $\Delta x_{i+1}$, $\Delta u_i$, and $\Delta \lambda_i$ are derived at each stage $i \in I_k$ for $\Delta x_i$, $\Delta t_{k}$, and $\Delta t_{k-1}$, respectively, from (\ref{eq:stateEquationLinearized})--(\ref{eq:LuLinearized}).
Moreover, variables are recursively eliminated from (\ref{eq:STOLinearized}) using these expressions.
When the stage of interest changes from $k$ to $k-1$, $\Delta t_k$ is further eliminated from (\ref{eq:STOLinearized}). That is, an expression of $\Delta t_k$ is derived with respect to $\Delta x_{i_k}$ and $\Delta t_{k-1}$, where $i_k := \min{\mathcal{I}_k}$.
This can be seen as the extension of the standard Riccati recursion that also recursively derives expressions of $\Delta x_{i+1}$, $\Delta u_i$, and $\Delta \lambda_i$ with respect to $\Delta x_i$ \cite{bib:mpcbook, bib:RiccatiMPC}.

\subsubsection{Terminal stage}
From (\ref{eq:KKTSystemTerminal}) at stage $N$, we obtain:
\begin{subequations}\label{subeq:riccati:terminal}
\begin{equation}\label{eq:lambdaN}
    \Delta \lambda_N = P_N \Delta x_N - s_N,
\end{equation}
where 
\begin{equation}
    P_N = Q_{xx, N}, \;\;
    s_N = - \bar{l}_{x, N}.
\end{equation}
\end{subequations}

\subsubsection{Intermediate stages}\label{subsubseq:intermediate}
Consider $i, i+1 \in \mathcal{I}_{k}$, $k \in \mathcal{K}$. 
Suppose that we have the expression of $\Delta \lambda_{i+1}$ with respect to $\Delta x_{i+1}$, $\Delta t_{k}$, and $\Delta t_{k-1}$ as
\begin{subequations}\label{subeq:riccati}
\begin{align}\label{eq:costate}
    \Delta \lambda_{i+1} = \; & P_{i+1} \Delta x_{i+1} - s_{i+1} + \Psi_{i+1} (\Delta t_k - \Delta t_{k-1}) \notag \\ 
    & + \Phi_{i+1} (- \Delta t_{k}),
\end{align}
where $P_{i+1} \in \mathbb{R}^{n_x \times n_x}$ and $s_{i+1}, \; \Psi_{i+1}, \; \Phi_{i+1} \in \mathbb{R}^{n_x}$.
Suppose also that we have the equations (\ref{eq:STOLinearized}) for $k$ and $k-1$ in which $\Delta x_{j}$ for $j \geq i+2$, $\Delta u_{j}$ for $j \geq i+1$, and $\Delta \lambda_{j}$ for $j \geq i+1$ are eliminated, that is, the equations (\ref{eq:STOLinearized}) for $k$ and $k-1$ are reduced to
\begin{align}\label{eq:STO1}
    & \sum_{j \in \mathcal{I}_k, \; j \leq i} \left( h_{x, j} ^{\rm T} \Delta x_j + h_{u, j} ^{\rm T} \Delta u_j + f_j ^{\rm T} \Delta \lambda_{j+1} + \bar{h}_j \right) \notag \\ 
    & + \Psi_{i+1} ^{\rm T} \Delta x_{i+1} + \xi_{i+1} (\Delta t_k - \Delta t_{k-1}) + \chi_{i+1} (- \Delta t_{k}) + \eta_{i+1} \notag \\ 
    & - \Phi_{i+1} ^{\rm T} \Delta x_{i+1} - \chi_{i+1} (\Delta t_k - \Delta t_{k-1}) - \rho_{i+1} (- \Delta t_{k}) - \iota_{i+1} \notag \\ 
   & = 0,
\end{align}
and
\begin{align}\label{eq:STO2}
    & \sum_{j \in \mathcal{I}_{k-1}} \left( h_{x, j} ^{\rm T} \Delta x_j + h_{u, j} ^{\rm T} \Delta u_j + f_j ^{\rm T} \Delta \lambda_{j+1} + \bar{h}_j \right) \notag \\ 
    & - \sum_{j \in \mathcal{I}_{K}, \; j \leq i} \left( h_{x, j} ^{\rm T} \Delta x_j + h_{u, j} ^{\rm T} \Delta u_j + f_j ^{\rm T} \Delta \lambda_{j+1} + \bar{h}_j \right) \notag \\ 
    & - \Psi_{i+1} ^{\rm T} \Delta x_{i+1} - \xi_{i+1} (\Delta t_k - \Delta t_{k-1}) - \chi_{i+1} (- \Delta t_{k}) \notag \\ 
    & - \eta_{i+1} = 0,
\end{align}
\end{subequations}
where $\xi_{i+1}, \chi_{i+1}, \rho_{i+1}, \eta_{i+1}, \iota_{i+1} \in \mathbb{R}$.
Note that $\Psi_{i+1} = 0$, $\Phi_{i+1} = 0$, $\xi_{i+1} = 0$, $\chi_{i+1} = 0$, $\rho_{i+1} = 0$, $\eta_{i+1} = 0$, and $\iota_{i+1} = 0$ at stage $i = N-1$.
Moreover, $\Psi_{i+1} = 0$, $\xi_{i+1} = 0$, $\chi_{i+1} = 0$, and $\eta_{i+1} = 0$ at stages $i = \min{\mathcal{I}_k} - 1$ and $k \in \left\{ 2, ..., K \right\}$, as explained in \ref{subsubsec:phaseTrans}.
First, the following equations are introduced:
\begin{subequations}\label{subeq:riccati:factorization}
\begin{equation}\label{eq:ricccati:F}
    F_i := Q_{xx, i} + A_i ^{\rm T} P_{i+1} A_i,
\end{equation}
\begin{equation}\label{eq:ricccati:H}
    H_i := Q_{xu, i} + A_i ^{\rm T} P_{i+1} B_i, 
\end{equation}
\begin{equation}\label{eq:ricccati:G}
    G_i := Q_{uu, i} + B_i ^{\rm T} P_{i+1} B_i,
\end{equation}
\begin{equation}
    \psi_{x, i} := h_{x, i} + A_i ^{\rm T} P_{i+1} f_i + A_i ^{\rm T} \Psi_{i+1}, 
\end{equation}
\begin{equation}
    \psi_{u, i} := h_{u, i} + B_i ^{\rm T} P_{i+1} f_s + B_i ^{\rm T} \Psi_{i+1},
\end{equation}
and
\begin{equation}
    \phi_{x, i} := A_i ^{\rm T} \Phi_{i+1}, \;\; 
    \phi_{u, i} := B_i ^{\rm T} \Phi_{i+1}.
\end{equation}
\end{subequations}
Second, $\Delta \lambda_{i+1}$ and $\Delta x_{i+1}$ are eliminated from (\ref{eq:LuLinearized}) using (\ref{eq:costate}) and (\ref{eq:stateEquationLinearized}). Subsequently, we can express $\Delta u_i$ using $\Delta x_{i}$, $\Delta t_k$, and $\Delta t_{k-1}$ as:
\begin{subequations}\label{subeq:riccati_u}
\begin{equation}\label{eq:riccati_u}
    \Delta u_i = K_i \Delta x_i + k_i + T_i (\Delta t_k - \Delta t_{k-1}) + W_i (- \Delta t_k), 
\end{equation}
where 
\begin{equation}
    K_i := - G_i ^{-1} H_i ^{\rm T},
\end{equation}
\begin{equation}
    k_i := - G_i ^{-1} (B_i ^{\rm T} P_{i+1} \bar{x}_i - B_i ^{\rm T} z_{i+1} + \bar{l}_{u, i}), 
\end{equation}
and
\begin{equation}
    T_{i} := - G_i ^{-1} \psi_{u, i}, \;\; W_i := - G_i ^{-1} \phi_{u, i}.
\end{equation}
\end{subequations}
Third, $\Delta \lambda_{i+1}$, $\Delta x_{i+1}$, and $\Delta u_i$ are eliminated from (\ref{eq:LxLinearized}) using (\ref{eq:costate}), (\ref{eq:stateEquationLinearized}), and (\ref{eq:riccati_u}), respectively. As a result, $\Delta \lambda_i$ using $\Delta x_{i}$, $\Delta t_k$, and $\Delta t_{k-1}$ is expressed as:
\begin{subequations}\label{subeq:riccati_lmd}
\begin{equation}\label{eq:riccati_lmd}
    \Delta \lambda_{i} = P_i \Delta x_i - s_i + \Psi_i (\Delta t_k - \Delta t_{k-1}) + \Phi_i (- \Delta t_{k}),
\end{equation}
where 
\begin{equation}
    P_i := F_i - K_i ^{\rm T} G_i K_i,
\end{equation}
\begin{equation}
    s_i := - \left\{ \bar{l}_{x, i} + A_i ^{\rm T} (P_{i+1} \bar{x}_i - s_{i+1}) + H_i k_i \right\},
\end{equation}
and 
\begin{equation}
    \Psi_i := \psi_{x, i} + K_i \psi_{u, i}, \;\;  \Phi_i := \phi_{x, i} + K_i \phi_{u, i}.
\end{equation}
\end{subequations}
Moreover, by eliminating $\Delta \lambda_{i+1}$, $\Delta x_{i+1}$, and $\Delta u_i$ from (\ref{eq:STO1}) and (\ref{eq:STO2}) using (\ref{eq:costate}), (\ref{eq:stateEquationLinearized}), and (\ref{eq:riccati_u}), we obtain: 
\begin{subequations}
\begin{align}\label{eq:STOres1}
    & \sum_{j \in \mathcal{I}_k, \; j \leq i-1} \left( h_{x, j} ^{\rm T} \Delta x_j + h_{u, j} ^{\rm T} \Delta u_j + f_j ^{\rm T} \Delta \lambda_{j+1} + \bar{h}_j \right) \notag \\ 
    & + \Psi_{i} ^{\rm T} \Delta x_{i} + \xi_{i} (\Delta t_k - \Delta t_{k-1}) + \chi_{i} (- \Delta t_{k}) + \eta_{i} \notag \\ 
    & - \Phi_{i} ^{\rm T} \Delta x_{i} - \chi_{i} (\Delta t_k - \Delta t_{k-1}) - \rho_{i} (- \Delta t_{k}) - \iota_{i} = 0
\end{align}
and
\begin{align}\label{eq:STOres2}
    & \sum_{j \in \mathcal{I}_{k-1}} \left( h_{x, j} ^{\rm T} \Delta x_j + h_{u, j} ^{\rm T} \Delta u_j + f_j ^{\rm T} \Delta \lambda_{j+1} + \bar{h}_j \right) \notag \\ 
    & - \sum_{j \in \mathcal{I}_{K}, \; j \leq i-1} \left( h_{x, j} ^{\rm T} \Delta x_j + h_{u, j} ^{\rm T} \Delta u_j + f_j ^{\rm T} \Delta \lambda_{j+1} + \bar{h}_j \right) \notag \\ 
    & - \Psi_{i} ^{\rm T} \Delta x_{i} - \xi_{i} (\Delta t_k - \Delta t_{k-1}) - \chi_{i} (- \Delta t_{k}) - \eta_{i} = 0,
\end{align}
\end{subequations}
where
\begin{subequations}\label{subeq:riccati_sto}
\begin{equation}
    \xi_i := \xi_{i+1} + f_i ^{\rm T} (P_{i+1} f_i + 2 \Psi_{i+1} ) + \psi_{u, i} ^{\rm T} T_i ,
\end{equation}
\begin{equation}
    \eta_i := \eta_{i+1} + \bar{h}_i + f_i ^{\rm T} (P_{i+1} \bar{x}_i - s_{i+1}) + \Psi_{i+1} ^{\rm T} \bar{x}_i + \psi_{u, i} ^{\rm T} k_i,
\end{equation}
\begin{equation}
    \chi_i := \chi_{i+1} + \Phi_{i+1} ^{\rm T} f_i + \psi_{u, i} ^{\rm T} W_i ,
\end{equation}
\begin{equation}
    \rho_i := \rho_{i+1} + \phi_{u, i} ^{\rm T} W_i ,
\end{equation}
and
\begin{equation}
    \iota_i := \iota_{i+1} + \Phi_{i+1} ^{\rm T} \bar{x}_i + \phi_{u, i} ^{\rm T} k_i .
\end{equation}
\end{subequations}
Therefore, we have equations (\ref{eq:riccati_lmd}), (\ref{eq:STOres1}), and (\ref{eq:STOres2}) for stage $i$ in the same form as equations (\ref{subeq:riccati}) for stage $i+1$.

\subsubsection{Phase transition stages}\label{subsubsec:phaseTrans}
Here, we consider phase $k \in \left\{ 1, ..., K \right\}$ and stage $i_k := \min \mathcal{I}_{k}$.
In this stage, the phase of interest changes from $k$ to $k-1$ when $k > 1$, and the backward recursion terminates when $k = 1$.
The equations (\ref{subeq:riccati}) until $i_k$ are expressed as:
\begin{subequations}
\begin{equation}\label{eq:costateTrans}
    \Delta \lambda_{i_{k}} = P_{i_{k}} \Delta x_{i_{k}} - s_{i_{k}} + \Psi_{i_{k}} (\Delta t_k - \Delta t_{k-1}) + \Phi_{i_{k}} (- \Delta t_k), 
\end{equation}
\begin{align}\label{eq:STO1Trans}
    & \Psi_{i_k} ^{\rm T} \Delta x_{i_k} + \xi_{i_k} (\Delta t_k - \Delta t_{k-1}) + \chi_{i_k} (- \Delta t_{k}) + \eta_{i_k} \notag \\ 
    & - \Phi_{i_k} ^{\rm T} \Delta x_{i_k} - \chi_{i_k} (\Delta t_k - \Delta t_{k-1}) - \rho_{i_k} (- \Delta t_{k}) - \iota_{i_k} = 0, 
\end{align}
and
\begin{align}\label{eq:STO2Trans}
    & \sum_{j \in \mathcal{I}_{k-1}} \left( h_{x, j} ^{\rm T} \Delta x_j + h_{u, j} ^{\rm T} \Delta u_j + f_j ^{\rm T} \Delta \lambda_{j+1} + \bar{h}_j \right) - \Psi_{i_k} ^{\rm T} \Delta x_{i_k} \notag \\ 
    & - \xi_{i_k} (\Delta t_k - \Delta t_{k-1}) - \chi_{i_k} (- \Delta t_{k}) - \eta_{i_k} = 0 .
\end{align}
\end{subequations}
Note that $\Delta t_{k-1} = \Delta t_0 = 0$ when $k=1$.
Therefore, from (\ref{eq:STO1Trans}), $\Delta t_{k}$ is determined as:
\begin{align}\label{eq:Deltatk}
    \Delta t_k = & - \sigma_{i_k} ^{-1} (\Psi_{i_k} - \Phi_{i_k}) ^{\rm T} \Delta x_{i_k} \notag \\ 
    & - \sigma_{i_k} ^{-1} (\xi_{i_k} - \chi_{i_k} ) (- \Delta t_{k-1})  - \sigma_{i_k} ^{-1} (\eta_{i_k} - \iota_{i_k} ).
\end{align}
Here, the following is defined: 
\begin{equation}\label{eq:sigma}
    \sigma_{i_k} := \xi_{i_k} - 2 \chi_{i_k} + \rho_{i_k} .
\end{equation}
The backward recursion is completed when $k=1$ because $i_1 = \min{\mathcal{I}_1} = 0$.
When $k > 1$, (\ref{eq:Deltatk}) is further substituted into (\ref{eq:costateTrans}) and (\ref{eq:STO2Trans}), which produces:  
\begin{subequations}\label{subeq:elimTk}
\begin{equation}\label{eq:afterElimTk_lmd}
    \Delta \lambda_{i_{K}} 
    = \tilde{P}_{i_{K}} \Delta x_{i_{K}} - \tilde{s}_{i_{K}} + \tilde{\Psi}_{i_{K}} (- \Delta t_{k-1})
\end{equation}
and
\begin{align}\label{eq:afterElimTk_STO}
    & \sum_{i \in \mathcal{I}_{K-1}} \left( h_{x, i} ^{\rm T} \Delta x_i + h_{u, i} ^{\rm T} \Delta u_i + f_i ^{\rm T} \Delta \lambda_{i+1} + \bar{h}_i \right) \notag \\ 
    & - \tilde{\Phi}_{i_k} ^{\rm T} \Delta x_{i_k} - \tilde{\rho}_{i_k} (- \Delta t_{k-1}) - \tilde{\iota}_{i_k} = 0,
\end{align}
where 
\begin{equation}\label{eq:afterElimTk_P}
    \tilde{P}_{i_k} := P_{i_k} - \sigma_{i_k} ^{-1} (\Psi_{i_k} - \Phi_{i_k}) (\Psi_{i_k} - \Phi_{i_k}) ^{\rm T} ,
\end{equation}
\begin{equation}\label{eq:afterElimTk_s}
    \tilde{s}_{i_k} := s_{i_k} + \sigma_{i_k} ^{-1} (\eta_{i_k} - \iota_{i_k}) (\Psi_{i_k} - \Phi_{i_k}) ,
\end{equation}
\begin{equation}\label{eq:afterElimTk_Phi}
    \tilde{\Phi}_{i_k} := \Psi_{i_k} - \sigma_{i_k} ^{-1} (\xi_{i_k} - \chi_{i_k}) (\Psi_{i_k} - \Phi_{i_k}) , 
\end{equation}
\begin{equation}\label{eq:afterElimTk_rho}
    \tilde{\rho}_{i_k} := \xi_{i_k} - \sigma_{i_k} ^{-1} (\xi_{i_k} - \chi_{i_k})^ 2 , 
\end{equation}
and
\begin{equation}\label{eq:afterElimTk_iota}
    \tilde{\iota}_{i_k} := \eta_{i_k} - \sigma_{i_k} ^{-1} (\xi_{i_k} - \chi_{i_k}) (\eta_{i_k} - \iota_{i_k}).
\end{equation}
\end{subequations}
Therefore, together with an equation (\ref{eq:STOLinearized}) for $k-1$, we obtain equations (\ref{eq:riccati_lmd}), (\ref{eq:STOres1}), and (\ref{eq:STOres2}) for stage $i_k$ in the same form as equations (\ref{subeq:riccati}) for stage $i+1$ with $P_{i+1} = \tilde{P}_{i_k}$, $s_{i+1} = \tilde{s}_{i_k}$, $\Psi_{i+1} = 0$, $\Phi_{i+1} = \tilde{\Psi}_{i_k}$, $\xi_{i+1} = 0$, $\chi_{i+1} = 0$, $\rho_{i+1} = \tilde{\rho}_{i_k}$, $\eta_{i+1} = 0$, and $\iota_{i+1} = \tilde{\iota}_{i_k}$.

\subsection{Forward Recursion}
After the backward recursion up to the initial stage ($i = 0$), all Newton directions are recursively computed from stage $0$ to $N$ forward in time using the results of the backward recursion.
First, the initial state direction $\Delta x_0$ is computed from (\ref{eq:initialStateLinearized}).
Second, the direction of the first switching time $\Delta t_1$ is computed from (\ref{eq:Deltatk}) with $\Delta t_{0} = 0$.
Third, $\Delta \lambda_i$, $\Delta u_i$, and $\Delta x_{i+1}$ are computed forward in time in each phase $k$ using (\ref{eq:costate}), (\ref{eq:riccati_u}), and (\ref{eq:stateEquationLinearized}), respectively.
Fourth, after this procedure, $\Delta t_{k+1}$ is computed using (\ref{eq:Deltatk}) until stage $i_{k+1} -1 = \min{\mathcal{I}_{k+1}} - 1$.
The third and forth steps are repeated from phase $k = 1$ to $k = K+1$, and the forward recursion is completed by computing $\Delta \lambda_{N}$ using (\ref{eq:lambdaN}).

\subsection{Properties of Proposed Riccati Recursion}
Next, we demonstrate that the proposed Riccati recursion always successfully solves the KKT system (\ref{subeq:KKTsystem}) if the solution is sufficiently close to the local minimum. This is in contrast to general QP solvers that fail to solve the KKT system because the Hessian matrix is inherently indefinite, as stated in Remark \ref{remark:indefiniteHessian}. 
We first explain that the proposed backward recursion is the same as dynamic programming for the QP subproblem (\ref{eq:QPsubproblem}), which is a quadratic approximation of the NLP (\ref{subeq:NLP}).
The following lemma states the equivalence of the dynamic programming and the proposed Riccati recursion:
\begin{lemma}\label{lemma:DPandRiccati}
Consider a phase $k \in \mathcal{K}$ and stage $i \in \mathcal{I}_k$. 
Suppose that $G_j \succ O$ for all $j \geq i$, and $\sigma_{i_l} > 0$ for all $l > k$, where $G_j$ and $\sigma_{i_l}$ are defined in (\ref{eq:ricccati:G}) and (\ref{eq:sigma}), respectively.
Then, the cost-to-go function of stage $i$ of dynamic programming for the QP (\ref{eq:QPsubproblem}) is expressed as: 
\begin{align}\label{eq:cost-to-go}
    \frac{1}{2} 
    & \begin{bmatrix}
        - \Delta t_k \\ 
        \Delta t_k - \Delta t_{k-1} \\ 
        \Delta x_{i} 
    \end{bmatrix}^{\rm T}
    \begin{bmatrix}
        \rho_{i} & \chi_{i} & \Phi_{i} ^{\rm T} \\
        \chi_{i} & \xi_{i} & \Psi_{i} ^{\rm T} \\
        \Phi_{i} & \Psi_{i} & P_{i}
    \end{bmatrix}
    \begin{bmatrix}
        - \Delta t_k \\ 
        \Delta t_k - \Delta t_{k-1} \\ 
        \Delta x_{i} 
    \end{bmatrix} \notag \\ 
    & + \begin{bmatrix}
        \iota_{i} \\ \eta_{i} \\ - s_{i} 
    \end{bmatrix}^{\rm T}
    \begin{bmatrix}
        - \Delta t_k \\ 
        \Delta t_k - \Delta t_{k-1} \\ 
        \Delta x_{i} 
    \end{bmatrix} ,
\end{align}
where $\xi_i$, $\rho_i$, $\Phi_i$, $\Psi_i$, $P_i$, $\iota_i$, $\eta_i$, and $s_i$ are defined by the Riccati recursion algorithm presented in subsection \ref{subsec:backwardRiccatiRecursion}.
\begin{figure*}
\begin{align}\label{eq:DPsubproblem}
    & \min_{\Delta u_i} \frac{1}{2} 
    \begin{bmatrix}
        - \Delta t_k \\ \Delta t_k - \Delta t_{k-1}  \\ \Delta x_i \\ \Delta u_i 
    \end{bmatrix}^{\rm T}
    \begin{bmatrix}
    \\
    & Q_{tt, i} & h_{x, i} ^{\rm T} & h_{u, i} ^{\rm T} \\
    & h_{x, i}  & Q_{xx, i} & Q_{xu, i} \\
    & h_{u, i}  & Q_{ux, i} & Q_{uu, i} \\
    \end{bmatrix}
    \begin{bmatrix}
        - \Delta t_k \\ \Delta t_k - \Delta t_{k-1}  \\ \Delta x_i \\ \Delta u_i 
    \end{bmatrix}
    + \begin{bmatrix}
        \\ \bar{h}_i \\ \bar{l}_{x, i} \\ \bar{l}_{u, i}
    \end{bmatrix}^{\rm T}
    \begin{bmatrix}
        - \Delta t_k \\ \Delta t_k - \Delta t_{k-1}  \\ \Delta x_i \\ \Delta u_i 
    \end{bmatrix} 
    \notag \\ 
    & \;\;\;\;\;\;\; + \frac{1}{2} 
    \begin{bmatrix}
        - \Delta t_k \\ 
        \Delta t_k - \Delta t_{k-1} \\ 
        \Delta x_{i+1} 
    \end{bmatrix}^{\rm T}
    \begin{bmatrix}
        \rho_{i+1} & \chi_{i+1} & \Phi_{i+1} ^{\rm T} \\
        \chi_{i+1} & \xi_{i+1} & \Psi_{i+1} ^{\rm T} \\
        \Phi_{i+1} &  \Psi_{i+1} & P_{i+1}
    \end{bmatrix}
    \begin{bmatrix}
        - \Delta t_k \\ 
        \Delta t_k - \Delta t_{k-1} \\ 
        \Delta x_{i+1} 
    \end{bmatrix}
    + \begin{bmatrix}
        \iota_{i+1} \\  \eta_{i+1} \\ - s_{i+1} 
    \end{bmatrix}^{\rm T}
    \begin{bmatrix}
        - \Delta t_k \\ 
        \Delta t_k - \Delta t_{k-1} \\ 
        \Delta x_{i+1} 
    \end{bmatrix} \notag \\
    & {\rm s.t.} \;\;\;\;\; A_i \Delta x_i + B_i \Delta u_i + f_i (\Delta t_k - \Delta t_{k-1}) - \Delta x_{i+1} + \bar{x}_i = 0,
\end{align}  
\end{figure*}
Moreover, the subproblem of the dynamic programming to determine $\Delta u_i$ is represented by (\ref{eq:DPsubproblem}), and determining $\Delta t_k$ is expressed as: 
\begin{equation}\label{eq:DPsubproblemTk}
    \min_{\Delta t_k} \;\; (\ref{eq:cost-to-go}) .
\end{equation}
\end{lemma}
\begin{proof}
The proof is achieved by induction. 
At the terminal stage ($i = N$), the cost-to-go function is represented by (\ref{eq:cost-to-go}) with $\Psi_N = \Phi_N = 0$ and $\xi_N = \chi_N = \rho_N = \eta_N = \iota_N = 0$.
Next, suppose that we have the cost-to-go function (\ref{eq:cost-to-go}) of stage $i+1$.
Then, it is clear that the subproblem of the dynamic programming is represented by (\ref{eq:DPsubproblem}). 
$\Delta u_i$ can be uniquely determined from (\ref{eq:DPsubproblem}) as (\ref{subeq:riccati_u}) because $G_i \succ O$. The cost-to-go function of stage $i$ is represented by (\ref{eq:cost-to-go}), where $P_i$, $s_i$, $\Psi_i$, $\Phi_i$, $\xi_i$, $\chi_i$, $\rho_i$, $\eta_i$, and $\iota_i$ are defined as (\ref{subeq:riccati:factorization}), (\ref{subeq:riccati_u}), (\ref{subeq:riccati_lmd}), and (\ref{subeq:riccati_sto}).
When $i = i_k := \min{I}_k$, $\Delta t_k$ can be further uniquely determined by solving (\ref{eq:DPsubproblemTk}) under the assumption $\sigma_{i_k} > 0$ after determining $\Delta u_{i_k}$ and obtaining the cost-to-go function of stage $i_k$ (\ref{eq:cost-to-go}). 
Then, the cost-to-go function of stage $i_k$ is in the form of (\ref{eq:cost-to-go}), where $P_{i_k}$, $s_{i_k}$, $\Phi_{i_k}$, $\rho_{i_k}$, and $\iota_{i_k}$ are defined as (\ref{subeq:elimTk}), $\Psi_{i_k} = 0$, $\xi_{i_k} = 0$, $\chi_{i_k} = 0$, and $\eta_{i_k} = 0$, respectively, which completes the proof.
\end{proof}
Note that the equivalence cannot be shown without the positive definiteness of $G_i$ and $\sigma_{i_k}$: if they are not positive definite, a unique solution does not exist and the cost-to-go function is not defined.
The following theorem is obtained based on Lemma \ref{lemma:DPandRiccati}.
Note that the discussion herein is restricted to the exact Hessian matrix to analyze the SOSC.
\begin{theorem}\label{theorem:exactSOSC}
We suppose that the SOSC and LICQ hold at the current iterate and 
consider that the exact Hessian matrix is used.
Then, $G_i$, as defined in (\ref{eq:ricccati:G}), is positive definite for all $i \in \left\{ 0, ..., N-1 \right\}$. 
Moreover, $\sigma_{i_k}$, as defined in (\ref{eq:sigma}), satisfies $\sigma_{i_k} > 0$ for all $k \in \left\{ 0, ..., K \right\}$.
\end{theorem}
\begin{proof}
The QP subproblem (\ref{eq:QPsubproblem}) with an exact Hessian matrix must have a unique global solution because the SOSC and LICQ hold \cite{bib:nocedal}. 
Therefore, the dynamic programming subproblems (\ref{eq:DPsubproblem}) and (\ref{eq:DPsubproblemTk}) must have unique solutions. 
Subsequently, the Hessian matrix with respect to $\Delta u_i$ in (\ref{eq:DPsubproblem}) must be positive definite, and the quadratic term with respect to $\Delta t_k$ in (\ref{eq:DPsubproblemTk}) must be positive.
Then, Lemma \ref{lemma:DPandRiccati} recursively shows that the Hessian matrix is $G_i$ with respect to $\Delta u_i$. Additionally, the lemma also shows that the quadratic term with respect to $\Delta t_k$ is $\sigma_{i_k}$, which completes the proof.
\end{proof}
Theorem \ref{theorem:exactSOSC} indicates that $G_i ^{-1}$ can always be efficiently computed using Cholesky factorizations if the current iterate is sufficiently close to a local minimum.
This fact also leads to the local convergence of the proposed method for NLP (\ref{subeq:NLP}) under the SOSC and LICQ. The proof for this is omitted because it is trivial.

\subsection{Reduced Hessian Modification via Riccati Recursion}
The recused Hessian matrix can be indefinite when the solution is not sufficiently close to a local minimum, such that the SQSC does not hold.
Subsequently, the KKT matrix is no longer invertible and the local convergence is not guaranteed. 
Efficient Cholesky factorization cannot be used to compute $G_i ^{-1}$.
Therefore, an algorithmic modification of the Riccati recursion is proposed to make the algorithm numerically robust and efficient for such situations, which can be considered a modification on the reduced Hessian matrix. 
To consider practical situations, Hessian approximations are allowed in the following, while Theorem \ref{theorem:exactSOSC} analyzes the exact Hessian matrix.
First, we introduce the following practical assumption:
\begin{assumption}\label{assump:HessianPositiveDefinite}
$\begin{bmatrix}
Q_{xx, i} & Q_{xu, i} \\ 
Q_{xu, i} ^{\rm T} & Q_{uu, i} 
\end{bmatrix} \succeq O$ and $Q_{uu, i} \succ O$ for all $i \in \left\{ 0, ..., N-1 \right\}$, $Q_{tt, k} \geq 0$ for all $k \in \left\{ 0, ..., K \right\}$, and $Q_{xx, N} \succeq O$.
\end{assumption}
This assumption is easily satisfied with the Gauss-Newton Hessian approximation or, more generally, with sequential convex programming  \cite{bib:SCPSurvey} for $\nabla_{x x} H_k (x_i, u_i, \lambda_{i+1})$, $\nabla_{x u} H_k (x_i, u_i, \lambda_{i+1})$, and $\nabla_{u u} H_k (x_i, u_i, \lambda_{i+1})$.
Under Assumption \ref{assump:HessianPositiveDefinite}, a unique solution of the dynamic programming subproblem (\ref{eq:DPsubproblem}) exists (that is, $G_i \succ O$) and $P_i \succeq O$ if $P_{i+1} \succeq O$, which is the same discussion as the dynamic programming for standard linear quadratic OCPs \cite{bib:DPandOCP}.
The solution to (\ref{eq:DPsubproblemTk}) also exists if $\sigma_{i_k} > 0$.
Based on these observations, a modification of the Riccati recursion algorithm is proposed, whereby $\tilde{P}_{i_k}$ is updated
at each phase-transition stage $i_k = \min{\mathcal{I}_k}$. This is expressed as: 
\begin{equation}\label{eq:Pmodification}
    \tilde{P}_{i_k} = P_{i_k}
\end{equation}
instead of (\ref{eq:afterElimTk_P}).
Then, $\tilde{P}_{i_k} \succeq O $, as long as ${P}_{i_k} \succeq O$.
Note that a different method of (\ref{eq:Pmodification}) can be used to make $\tilde{P}_{i_k} \succeq O$, which eliminates the negative curvature from $P_{i_k}$ through an eigenvalue decomposition \cite{bib:exactHessianReg}. However, this requires much more computational time than (\ref{eq:Pmodification}).
Nevertheless, the computationally cheap modification (\ref{eq:Pmodification}) works surprisingly well in practice.
We also modify $\sigma_{i_k}$ instead of (\ref{eq:sigma}), which is expressed as:
\begin{equation}\label{eq:sigmaModification}
    \tilde{\sigma}_{i_k} = \begin{cases}
        \sigma_{i_k} & \sigma_{i_k} > \sigma_{\rm \min} \\ 
        |\sigma_{i_k}| + \bar{\sigma} & \sigma_{i_k} \leq \sigma_{\rm \min} \\ 
    \end{cases},
\end{equation}
where $\sigma_{\rm \min} \geq 0$ and $\bar{\sigma} > 0$.
A practical rule to choose $\sigma_{\rm min}$ and $\bar{\sigma}$ is as follows.
We empirically observed that numerical ill-conditioning in the Newton-type method produced a large $\Delta t_s$, that is, a too small $\sigma_{i_k}$ in (\ref{eq:Deltatk}).
From this observation, a desired maximum absolute value of $\Delta t_s$ was chosen as $\Delta t_{s, {\rm max}} > 0$. 
Subsequently, $\sigma_{\rm min}$ and $\bar{\sigma}$ were chosen such that the absolute value of $\sigma_{i_k} ^{-1} (\eta_{i_k} - \iota_{i_k})$ did not exceed $\Delta t_{s, {\rm max}}$. This is expressed as: 
\begin{equation}\label{eq:sigmaModification2}
        \sigma_{\rm \min} = \bar{\sigma} = \Delta t_{s, {\rm max}} ^{-1} |\eta_{i_k} - \iota_{i_k}| 
\end{equation}
for each $k \in \mathcal{K}$.
For example, it was discovered that choosing $\Delta t_{s, {\rm max}}$ from a range of 0.1 to 1.0 worked well in practice.

The following Theorem claims that the proposed algorithmic modification makes the reduced Hessian matrix of the KKT system (\ref{subeq:KKTsystem}) positive definite and the KKT matrix invertible:
\begin{theorem}
Suppose that Assumption \ref{assump:HessianPositiveDefinite} holds.
Then, the reduced Hessian matrix of the KKT system (\ref{subeq:KKTsystem}) with the modifications (\ref{eq:Pmodification}) and (\ref{eq:sigmaModification}) is positive definite.
If, in addition, the LICQ holds, then the KKT matrix is invertible.
\end{theorem}
\begin{proof}
Under Assumption \ref{assump:HessianPositiveDefinite} and with the Hessian modifications (\ref{eq:Pmodification}) and (\ref{eq:sigmaModification}), $P_i \succeq O$ for all $i = 0, ..., N$, $G_i \succ O$ for all $i = 0, ..., N-1$, and $\sigma_{i_k} > 0$ for all $k = 1,..., K$ hold.
Therefore, the QP subproblems of the dynamic programming (\ref{eq:DPsubproblem}) and (\ref{eq:DPsubproblemTk}) always have unique global solutions. 
Theorem 16.2 of \cite{bib:nocedal} completes the first claim.
The second claim then directly follows from Theorem 16.1 of \cite{bib:nocedal}.
\end{proof}
\begin{remark}
The KKT matrix is always invertible under Assumption \ref{assump:HessianPositiveDefinite} and the proposed reduced Hessian modifications. Therefore, the local convergence of the Newton-type method is ensured if the resultant reduced Hessian is sufficiently close to the exact one \cite{bib:nocedal, bib:SCPSurvey} and the LICQ holds.
\end{remark}

\subsection{Algorithm}
The single Newton iteration using the proposed Riccati recursion algorithm is summarized in Algorithm \ref{alg2}.
After computing the KKT system (\ref{subeq:KKTsystem}) (line 1), the proposed method recursively eliminates the Newton steps from the KKT system (\ref{subeq:KKTsystem}) in the backward recursion (lines 3--14) and recursively computes the Newton steps in the forward recursion (lines 16--27).
Finally, the Newton steps of the slack variables and Lagrange multipliers are computed from the other Newton steps (line 28) according to the primal-dual interior point method  \cite{bib:nocedal, bib:ipopt}.
Calculations in the proposed algorithm, such as matrix inversions or matrix-matrix multiplications, do not grow with respect to the length of the horizon $N$. Therefore, the computational time of the proposed method is $O(N)$ as the standard Riccati recursion algorithm \cite{bib:RiccatiMPC}.

\begin{algorithm}[tb]
\caption{Computation of Newton step via proposed Riccati recursion}
\label{alg2}
\begin{algorithmic}[1]
    \Require Initial state ${x} (t_0)$ and the current iterate $x_0$, ..., $x_{N}$, $x_s$, $u_0$, ..., $u_{N-1}$, $u_s$, $\lambda_0$, ..., $\lambda_N$, $\lambda_s$, and $t_s$.
    \Ensure Newton directions $\Delta x_0$, ..., $\Delta x_{N}$, $\Delta x_s$, $\Delta u_0$, ..., $\Delta u_{N-1}$, $\Delta u_s$, $\Delta \lambda_0$, ..., $\Delta \lambda_N$, $\Delta \lambda_s$, and $\Delta t_s$.
    \State Compute the KKT system (\ref{subeq:KKTsystem}).
    \State \texttt{// Backward recursion}
    \State Compute $P_N$ and $z_N$ from (\ref{subeq:riccati:terminal}). 
    \For{$k=K+1, \cdots, 1$} 
        \For{$i=\max{ \mathcal{I}_k}, \cdots, \min{ \mathcal{I}_k }$} 
            \State Compute $K_i$, $k_i$, $T_i$, and $W_i$ from (\ref{subeq:riccati_u}) .
            \State Compute $P_i$, $z_i$, $\Psi_i$, and $\Phi_i$ from (\ref{subeq:riccati_lmd}). 
            \State Compute $\xi_i$, $\xi_i$, $\rho_i$, $\eta_i$, and $\iota_i$, from (\ref{subeq:riccati_sto}). 
        \EndFor
        \If{$k < K+1$}
            \State Compute $\tilde{P}_{i_k}$, $\tilde{s}_{i_k}$, $\tilde{\Phi}_{i_k}$, $\tilde{\rho}_{i_k}$, and $\tilde{\iota}_{i_k}$ from (\ref{subeq:elimTk}) optionally with modifications (\ref{eq:Pmodification}) and (\ref{eq:sigmaModification}). 
        \EndIf
        \State Set $\Psi_{i_k}$, $\xi_{i_k}$, $\chi_{i_k}$, and $\iota_{i_k}$ to zero. 
    \EndFor
    \State \texttt{// Forward recursion}
    \State Compute $\Delta x_0$ from (\ref{eq:initialStateLinearized}).
    \State Compute $\Delta t_1$ from (\ref{eq:Deltatk}).
    \For{$k=1, \cdots, K+1$} 
        \For{$i=\min{ \mathcal{I}_k}, \cdots, \max{ \mathcal{I}_k }$} 
            \State Compute $\Delta u_i$ and $\Delta \lambda_i$ from (\ref{subeq:riccati_u}) and (\ref{eq:riccati_lmd}), respectively.
            \State Compute $\Delta x_{i+1}$ from (\ref{eq:stateEquationLinearized})
        \EndFor
        \If{$k < K+1$}
            \State Compute $\Delta t_{k+1}$ from (\ref{eq:Deltatk})
        \EndIf
    \EndFor
    \State Compute $\Delta \lambda_N$ from (\ref{subeq:riccati:terminal}).
    \State Compute $\Delta z_i$, $\Delta \nu_i$, $\Delta w_i$, and $\Delta \upsilon_i$ for $i= 0, ..., N-1$ from the other Newton steps.
\end{algorithmic}
\end{algorithm}

\section{State Jumps and Switching Conditions}\label{sec:stateJumpsAndSwitchingConditions}
In this section, the proposed NLP formulations and Riccati recursion algorithm are further extended to switched systems involving state jumps and switching conditions, which have not yet been considered.
Some classes of switched systems involve state jumps at the same time as the switch, which is expressed as: 
\begin{equation}\label{eq:stateJump}
    x (t_{k}) = f_j (x (t_{k} -)), \;\; f_j : \mathbb{R}^{n_x} \to \mathbb{R}^{n_x},
\end{equation}
where $t_k-$ is the instant immediately before $t_k$.
We also consider that there is a stage cost on the state immediately before the state jump. That is, it is assumed that $l_j (x (t_k -))$ is added to the cost function (\ref{eq:costFunction}).
Moreover, such switches are often state-dependent, that is, the switch occurs if the state satisfies some condition (hereafter called the switching condition), which is expressed as:
\begin{equation}\label{eq:switchingCondition}
    e (x (t_k -)) = 0, \;\; e : \mathbb{R}^{n_x} \to \mathbb{R}^{n_e}.
\end{equation}
For example, a legged robot can be modeled as a system with switches involving state jumps and switching conditions. When the distance between the robot’s foot and the ground becomes zero (the switching condition), the generalized velocity instantly changes due to the impact (the state jump) \cite{bib:twoStageApplication1, bib:twoStageApplication2}.

If there is a state jump at the $k$th switch, 
a new grid point $i_k -$ is introduced corresponding to time $t_{k} -$, which is added to the KKT conditions and expressed as:  
\begin{subequations}
\begin{equation}\label{eq:NLP:jumpStateEquation}
    x_{i_k} = f_j (x_{i_k -} ) 
\end{equation}
and
\begin{equation}\label{eq:NLP:jumpCostateEquation}
    \nabla_x l_{j} (x_{i_k - } ) + \nabla_x f_j ^{\rm T} (x_{i_k - }) \lambda_{i_k} - \lambda_{i_k -} = 0. 
\end{equation}
\end{subequations}
The KKT systems regarding the state jump are therefore expressed as:
\begin{subequations}
\begin{equation}\label{eq:JumpEquationLinearized}
    \Delta x_{i_k} = A_{i_k -} \Delta x_{i_k -} + \bar{x}_{i_k -}
\end{equation}
and
\begin{equation}\label{eq:CostateJumpLienarized}
    \Delta \lambda_{i_k -} = Q_{xx, i_k -} \Delta x_{i_k -} + A_{i_k -} ^{\rm T} \Delta \lambda_{i_k} + \bar{l}_{x, i_k -} ,
\end{equation}
where $A_{i_k -} := \nabla_x f_j (x_{i_k -})$ and $Q_{xx, i_k -}$ is the Hessian of the Lagrangian with respect to $x_{i_k -}$. $\bar{x}_{i_k}$ and $\bar{l}_{x, i_k -}$ are the residuals in (\ref{eq:NLP:jumpStateEquation}) and (\ref{eq:NLP:jumpCostateEquation}), respectively.
\end{subequations}
When (\ref{eq:JumpEquationLinearized}) and (\ref{eq:CostateJumpLienarized}) are considered in the backward recursion, we have the equations until stage $i_k = \min{\mathcal{I}_k}$.
That is, the equations are in the form of (\ref{eq:afterElimTk_lmd}), (\ref{eq:afterElimTk_STO}), and (\ref{eq:STOLinearized}), with the phase index $k$ replaced with $k-1$.
By eliminating $\Delta \lambda_{i_k}$ and $\Delta x_{i_k}$ from these equations using (\ref{eq:JumpEquationLinearized}) and (\ref{eq:CostateJumpLienarized}), we then obtain: 
\begin{subequations}
\begin{equation}
    \Delta \lambda_{i_k -} = P_{i_k -} \Delta x_{i_k -} + \Phi_{i_k -} (- \Delta t_{k-1}) - s_{i_k -} , 
\end{equation}
\begin{align}
    & \sum_{i \in \mathcal{I}_{k-1}} \left( h_{x, i} ^{\rm T} \Delta x_i + h_{u, i} ^{\rm T} \Delta u_i + f_i ^{\rm T} \Delta \lambda_{i+1} + \bar{h}_i \right) \notag \\ 
    & - \Phi_{i_k -} ^{\rm T} \Delta x_{i_k -} - \xi_{i_k -} (\Delta t_{k-1}) - \eta_{i_k -} = 0,
\end{align}
and
\begin{align}
    & \sum_{i \in \mathcal{I}_{k-2}} \left( h_{x, i} ^{\rm T} \Delta x_i + h_{u, i} ^{\rm T} \Delta u_i + f_i ^{\rm T} \Delta \lambda_{i+1} + \bar{h}_i \right) \notag \\ 
    & - \sum_{i \in \mathcal{I}_{k-1}} \left( h_{x, i} ^{\rm T} \Delta x_i + h_{u, i} ^{\rm T} \Delta u_i + f_i ^{\rm T} \Delta \lambda_{i+1} + \bar{h}_i \right) = 0,
\end{align}
where 
\begin{equation}
    P_{i_k - } := Q_{xx, i_k -} + A_{i_k - } ^{\rm T} \tilde{P}_{i_k} A_{i_k -},
\end{equation}
\begin{equation}
    \Phi_{i_k - } := A_{i_k - } ^{\rm T} \tilde{\Phi}_{i_k},
\end{equation}
\begin{equation}
    s_{i_k -} := A_{i_k -} ^{\rm T} ( \tilde{s}_{i_k} - \tilde{P}_{i_k} \bar{x}_{i_k -}) - \bar{l}_{x, i_k -} ,
\end{equation}
and
\begin{equation}
    \xi_{i_k -} = \tilde{\xi}_{i_k}, \;\;\;
    \eta_{i_k - } = \tilde{\eta}_{i_k} + \tilde{\Phi}_{i_k} ^{\rm T} \bar{x}_{i_k -}.
\end{equation}
\end{subequations}
Therefore, the same equation as (\ref{subeq:riccati_lmd}) is achieved at stage $i_k-$, and we can proceed to stage $i_k-1$.
In the forward recursion, $\Delta x_{i_k -}$ is computed from $\Delta x_{i_k -1}$ and $\Delta u_{i_k -1}$ using (\ref{eq:stateEquationLinearized}) at stage $i_k - 1$. 
Then, $\Delta x_{i_k}$ and $\Delta \lambda_{i_k -}$ are computed from $\Delta x_{i_k - }$ according to (\ref{eq:JumpEquationLinearized}) and (\ref{eq:riccati_lmd}) for $i = i_k - $ with $\Psi_{i_k} = 0$.

\subsection{Switching Conditions}\label{subseq:switchingConditions}
The switching conditions are typically presented as pure-state equality constraints, which are difficult to treat efficiently with the Riccati recursion algorithm (specifically with the $O(N)$ computational burden) \cite{bib:constrainedRiccati}.
In this section, we assume that the state is partitioned into the generalized coordinate and velocity as $x = \begin{bmatrix}
        q ^{\rm T} &
        v ^{\rm T}
\end{bmatrix} ^{\rm T}$, $q, v \in \mathbb{R}^{n}$ and the state equation of each subsystem is expressed as: 
\begin{equation}\label{eq:stateEquationForm}
    f (x, u) := 
    \begin{bmatrix}
        f ^{(q)} (x)  \\ 
        f ^{(v)} (x, u) 
    \end{bmatrix}, 
\end{equation}
where $f ^{(q)} : \mathbb{R}^{n_x} \to \mathbb{R}^{n}$ and $f ^{(v)} : \mathbb{R}^{n_x} \times \mathbb{R}^{n_u} \to \mathbb{R}^{n}$. 
Moreover, we assume that the pure-state equality constraints representing the switching conditions have a relative degree of two, that is, the Jacobian of (\ref{eq:switchingCondition}) is expressed as: 
\begin{equation}\label{eq:pureStateConstraintJacobian}
    \nabla_x e (x) = 
    \begin{bmatrix}
        \nabla_q e (q) & O
    \end{bmatrix}, \;\; \nabla_q e (q) \in \mathbb{R}^{n_e \times n}.
\end{equation}
These assumptions mainly represent the position-level constraints on mechanical systems, such as robotic systems with rigid contacts.
To perform the Newton-type method with $O(N)$ complexity under these assumptions, we did not consider the pure-state equality constraint 
\begin{equation}
    e (x_{i_k} -)  = 0 
\end{equation}
directly as \cite{bib:constrainedRiccati}.
Instead, the constraint was transformed into a mixed state-input constraint at the two-stage grid point before the switch \cite{bib:stateConstrainedRiccati}, which is expressed as:
\begin{equation}\label{eq:transformedSwitchingCondition}
    e_{i} := e \left( x_{i} + f (x_{i}, u_{i}) \Delta \tau_k + g (x_{i} + f (x_{i}, u_{i}) \Delta \tau_k) \Delta \tau_k \right) = 0,
\end{equation}
where $i$ is two-stage before the grid $i_k -$, that is, $i$ satisfies $i + 2 = i_k -$.
The Lagrange multiplier is introduced with respect to (\ref{eq:transformedSwitchingCondition}), $\zeta_i \in \mathbb{R}^{n_e}$. 
The KKT conditions regarding this stage are then represented by (\ref{eq:NLP:stateEquation}), (\ref{eq:transformedSwitchingCondition}),
$\tilde{r}_{x, i} := r_{x, i} + C_i ^{\rm T} \zeta_i = 0$,
and $\tilde{r}_{u, i} := r_{u, i} + D_i ^{\rm T} \zeta_i = 0$,
where $C_i := \nabla_{x} e(\cdot)$ and $D_i := \nabla_{u} e(\cdot)$.
Furthermore, the KKT conditions of (\ref{eq:KKTSTO}) related to $t_{k-1}$ are replaced with:
\begin{align*}
    & E_i ^{\rm T} \zeta_i + \frac{1}{N_{k-1}} \sum_{j \in \mathcal{I}_{k-1}} H_{k-1} (x_j, u_j, \lambda_{j+1}) \notag \\ 
    & - \frac{1}{N_{k}} \sum_{i \in \mathcal{I}_{k}} H_{k} (x_j, u_j, \lambda_{j+1}) + \upsilon_{k-1} - \upsilon_{k} = 0 
\end{align*}
and
\begin{align*}
    & \frac{1}{N_{k-2}} \sum_{j \in \mathcal{I}_{k-2}} H_{k-2} (x_j, u_j, \lambda_{j+1}) - E_i ^{\rm T} \zeta_i \notag \\ 
    & - \frac{1}{N_{k-1}} \sum_{i \in \mathcal{I}_{k-1}} H_{k-1} (x_j, u_j, \lambda_{j+1}) + \upsilon_{k-2} - \upsilon_{k-1} = 0 ,
\end{align*}
where $E_i := \frac{1}{N_{k-1}} \nabla_{t_k} e(\cdot)$.
The KKT systems regarding this stage are then represented by (\ref{eq:stateEquationLinearized}) and expressed as: 
\begin{subequations}
\begin{equation}
    C_{i} \Delta x_{i} + D_{i} \Delta u_{i} + E_{i} (\Delta t_k - \Delta t_{k-1}) + \bar{e}_{i} = 0, 
\end{equation}
\begin{align}
    & Q_{xx, i} \Delta x_i + Q_{xu, i} \Delta u_i + A_i ^{\rm T} \Delta \lambda_{i+1} + C_i ^{\rm T} \Delta \zeta_i - \Delta \lambda_{i} \notag \\ 
    & + h_{x, i} (\Delta t_k - \Delta t_{k-1}) + \tilde{l}_{x, i} = 0,
\end{align}
and
\begin{align}
    & Q_{xu, i} ^{\rm T} \Delta x_i + Q_{uu, i} \Delta u_i + B_i ^{\rm T} \Delta \lambda_{i+1} + D_i ^{\rm T} \Delta \zeta_i \notag \\ 
    & + h_{u, i} (\Delta t_k - \Delta t_{k-1}) + \tilde{l}_{u, i} = 0,
\end{align}
where 
\begin{equation*}
    \tilde{l}_{x, i} := \tilde{r}_{x, i} + \nabla_{x} g ^{\rm T} (x_i, u_i) {\rm diag}(z_i) ^{-1} ({\rm diag}(\nu_i) r_{g, i} - r_{z, i})
\end{equation*}
and
\begin{equation*}
    \tilde{l}_{u, i} := \tilde{r}_{u, i} + \nabla_{u} g ^{\rm T} (x_i, u_i) {\rm diag}(z_i) ^{-1} ({\rm diag}(\nu_i) r_{g, i} - r_{z, i}).
\end{equation*}
In addition, $\bar{e}_i$ is residual in (\ref{eq:transformedSwitchingCondition}). 
\end{subequations}
In the backward recursion at stage $i$, we have (\ref{eq:costate}),
\begin{subequations}\label{subeq:switchingConditionSTO}
\begin{align}\label{eq:switchingConditionSTOeq1}
    & \sum_{j \in \mathcal{I}_{k-1}, \; j \leq i} \left( h_{x, j} ^{\rm T} \Delta x_j + h_{u, j} ^{\rm T} \Delta u_j + f_j ^{\rm T} \Delta \lambda_{j+1} + \bar{h}_j \right) \notag \\ 
    & + E_i ^{\rm T} (\Delta \zeta_i + \zeta_i) + \Psi_{i+1} ^{\rm T} \Delta x_{i+1} + \xi_{i+1} (\Delta t_{k-1} - \Delta t_{k-2}) \notag \\ 
    & + \chi_{i+1} (- \Delta t_{k-1}) + \eta_{i+1} - \Phi_{i+1} ^{\rm T} \Delta x_{i+1} \notag \\ 
    & - \chi_{i+1} (\Delta t_{k-1} - \Delta t_{k-2}) - \rho_{i+1} (- \Delta t_{k-1}) - \iota_{i+1} = 0,
\end{align}
and 
\begin{align}\label{eq:switchingConditionSTOeq2}
    & \sum_{j \in \mathcal{I}_{k-2}} \left( h_{x, j} ^{\rm T} \Delta x_j + h_{u, j} ^{\rm T} \Delta u_j + f_j ^{\rm T} \Delta \lambda_{j+1} + \bar{h}_j \right) \notag \\ 
    & - \sum_{j \in \mathcal{I}_{K-1}, \; j \leq i} \left( h_{x, j} ^{\rm T} \Delta x_j + h_{u, j} ^{\rm T} \Delta u_j + f_j ^{\rm T} \Delta \lambda_{j+1} + \bar{h}_j \right) \notag \\ 
    & - \Psi_{i+1} ^{\rm T} \Delta x_{i+1} - \xi_{i+1} (\Delta t_{k-1} - \Delta t_{k-2}) - \chi_{i+1} (- \Delta t_{k-1}) \notag \\ 
    & - \eta_{i+1} + E_i ^{\rm T} (\Delta \zeta_i + \zeta_i) = 0.
\end{align}
\end{subequations}
$\Delta u_i$ and $\Delta \zeta_i$ are then eliminated, resulting in: 
\begin{subequations}
\begin{align}\label{eq:DeltauAndDeltanu}
    \begin{bmatrix} \Delta u_{i} \\ \Delta \zeta_i \end{bmatrix} = \; & \begin{bmatrix} K_{i} ^{\rm T} \\ M_i \end{bmatrix} \Delta x_{i} 
    + \begin{bmatrix} T_{i} \\ L_i \end{bmatrix} (\Delta t_{k-1} - \Delta t_{k-2}) \notag \\ 
    & + \begin{bmatrix} W_{i} \\ N_i \end{bmatrix} (- \Delta t_{k-1}) 
    + \begin{bmatrix} k_{i} \\ m_i \end{bmatrix} ,
\end{align}
where 
\begin{equation}
    \begin{bmatrix} K_{i} \\ M_i \end{bmatrix}
    := - \begin{bmatrix} G_{i} & D_i ^{\rm T} \\ D_i & \end{bmatrix}^{-1} 
        \begin{bmatrix} H_{i} ^{\rm T} \\ C_i \end{bmatrix} ,
\end{equation}
\begin{equation}
    \begin{bmatrix} T_{i} \\ L_i \end{bmatrix}
    := - \begin{bmatrix} G_{i} & D_i ^{\rm T} \\ D_i & \end{bmatrix}^{-1} 
        \begin{bmatrix} \psi_{u, i} \\ E_i \end{bmatrix} ,
\end{equation}
\begin{equation}
    \begin{bmatrix} W_{i} \\ N_i \end{bmatrix}
    := - \begin{bmatrix} G_{i} & D_i ^{\rm T} \\ D_i & \end{bmatrix}^{-1} 
        \begin{bmatrix} \phi_{u, i} \\ 0 \end{bmatrix} , 
\end{equation}
and
\begin{equation}
    \begin{bmatrix} k_{i} \\ m_i \end{bmatrix}
    := - \begin{bmatrix} G_{i} & D_i ^{\rm T} \\ D_i & \end{bmatrix}^{-1} 
        \begin{bmatrix} {B_{i}}^{\rm T} (P_{i+1} \bar{x}_{i} - s_i) + \tilde{l}_{u, i} \\ \bar{e}_i \end{bmatrix} .
\end{equation}
\end{subequations}
Therefore, we obtain the form equations of (\ref{eq:riccati_lmd}), (\ref{eq:STOres1}), and (\ref{eq:STOres2}) for stage $i$ with the following:
\begin{subequations}
\begin{equation}
    P_i = F_i - \begin{bmatrix}
        K_i \\ M_i 
    \end{bmatrix}^{\rm T}
    \begin{bmatrix} G_{i} & D_i ^{\rm T} \\ D_i & \end{bmatrix}
    \begin{bmatrix}
        K_i \\ M_i 
    \end{bmatrix},
\end{equation}
\begin{equation}
    s_i = - \left\{ \tilde{l}_{x, i} + A_i ^{\rm T} (P_{i+1} \bar{x}_i - s_{i+1}) + H_i k_i + C_i ^{\rm T} m_i \right\},
\end{equation}
and
\begin{equation}
    \Psi_i = \psi_{x, i} 
    +  \begin{bmatrix}
        K_i \\ M_i 
    \end{bmatrix}^{\rm T}
    \begin{bmatrix}
        \psi_{u, i} \\ E_i 
    \end{bmatrix}, \;\;
    \Phi_i = \phi_{x, i} + \begin{bmatrix}
        K_i \\ M_i 
    \end{bmatrix}^{\rm T}
    \begin{bmatrix}
        \phi_{u, i} \\ 0
    \end{bmatrix}.
\end{equation}
\end{subequations}
\begin{subequations}
\begin{equation}
    \xi_i = \xi_{i+1} + f_i ^{\rm T} (P_{i+1} f_i + 2 \Psi_{i+1} ) + \psi_{u, i} ^{\rm T} T_i + E_i ^{\rm T} L_i,
\end{equation}
\begin{align}
    \eta_i = \; & \eta_{i+1} + \bar{h}_i + E_i ^{\rm T} \zeta_i + f_i ^{\rm T} (P_{i+1} \bar{x}_i - s_{i+1}) + \Psi_{i+1} ^{\rm T} \bar{x}_i \notag \\ 
    & + \psi_{u, i} ^{\rm T} k_i + E_i ^{\rm T} m_i,
\end{align}
\begin{equation}
    \chi_i = \chi_{i+1} + \Phi_{i+1} ^{\rm T} f_i + \psi_{u, i} ^{\rm T} W_i + E_i ^{\rm T} N_i,
\end{equation}
\begin{equation}
    \rho_i = \rho_{i+1} + \phi_{u, i} ^{\rm T} W_i ,
\end{equation}
and
\begin{equation}
    \iota_i = \iota_{i+1} + \Phi_{i+1} ^{\rm T} \bar{x}_i + \phi_{u, i} ^{\rm T} k_i.
\end{equation}
\end{subequations}
At the forward recursion, $\Delta u_i$ and $\Delta \zeta_i$ are computed using (\ref{eq:DeltauAndDeltanu}) instead of only computing $\Delta u_i$ as (\ref{eq:riccati_u}).

\section{Numerical Experiments}\label{sec:experiments}
\subsection{Comparison with Off-The-Shelf Solvers}
\subsubsection{Problem settings}
The effectiveness of the proposed method was demonstrated through two numerical experiments.
The first experiment was a comparison with off-the-shelf NLP solvers on the switched system, consisting of three nonlinear subsystems treated in \cite{bib:two-stage0, bib:twoStage}. The dynamics of the subsystems are expressed as: 
\begin{equation*}
    f_1 (x, u) = 
    \begin{bmatrix}
    x_1 + u_1 \sin(x_1) \\
    - x_2 - u_1 \cos(x_2) 
    \end{bmatrix}, 
\end{equation*}
\begin{equation*}
    f_2 (x, u) = 
    \begin{bmatrix}
    x_2 + u_1 \sin(x_2) \\
    - x_1 - u_1 \cos(x_1) 
    \end{bmatrix}, 
\end{equation*}
and
\begin{equation*}
    f_3 (x, u) = 
    \begin{bmatrix}
    - x_1 - u_1 \sin(x_1) \\
    x_2 + u_1 \cos(x_2) 
    \end{bmatrix}.
\end{equation*}
The stage cost is expressed as:
\begin{equation*}
    l_{k} (x, u) = \frac{1}{2} ||x - x_{\rm ref} ||^2 + || u ||^2 , \;\; k \in \left\{ 1, 2 ,3 \right\}.
\end{equation*}
In addition, the terminal cost is expressed as:
\begin{equation*}
    V_{f} (x) = \frac{1}{2} ||x - x_{\rm ref} ||^2 , 
\end{equation*}
where $x_{\rm ref} = [1, \; -1]^{\rm T}$.
The initial and terminal times of the horizon are denoted by $t_0 = 0$ and $t_3 = 3$, respectively.
The initial state is denoted by $x (t_0) = [2, \; 3] ^{\rm T}$.

Benchmarks were set for the solution times of the NLP (\ref{subeq:NLP}) example against the {Ipopt} \cite{bib:ipopt} and a SQP methods with {qpOASES} \cite{bib:qpOASES}, both of which were used through the {CasADi} \cite{bib:casadi} interface.
The mesh-refinement was not considered in this comparison in order to focus on a pure comparison of the abilities to solve the NLP problems. 
That is, only the single NLP step of Algorithm \ref{alg1} was considered.
{Ipopt} was used with the default settings of {CasADi}. For example, MUMPS \cite{bib:MUMPS1} was used to solve the linear systems.
The built-in SQP solver of {CasADi} was used with the exact Hessian matrix and qpOASES as the backend QP solver.
The Hessian regularization was enabled within the {qpOASES} options and the Hessian type was set to ``indefinite".
The proposed method was written in C++, used Eigen \cite{bib:EigenWeb} for the linear algebra, and used an exact Hessian matrix and reduced Hessian modifications (\ref{eq:Pmodification}) and (\ref{eq:sigmaModification}) throughout the experiments. 
$\sigma_{\rm min}$ and $\bar{\sigma}$ were chosen using (\ref{eq:sigmaModification2}) with $\Delta t_{k, {\rm max}} = 0.5$.
The proposed method only used the fraction-to-boundary rule \cite{bib:nocedal, bib:ipopt} for the step size selection and monotonically decreased the barrier parameter $\epsilon$ when the $l_2$ norm of the KKT residual (in the perturbed KKT conditions (\ref{subeq:KKTconditions})) was smaller than 0.1.
The solution times were compared for these algorithms with the total number of grids set as $N = 10$, $50$, $100$, and $500$. 
$N$ was divided into $N_1$, $N_2$, and $N_3$, so that each phase had an almost equal number of discretization grids. The values of $[N_1, N_2, N_3]$ were set to $[4, 3, 3]$,  $[17, 17, 16]$, $[34, 33, 33]$, and $[167, 167, 166]$ for $N = 10$, $50$, $100$, and $500$, respectively. 
The minimum dwell-times were set as $\underline{\Delta}_1 = \underline{\Delta}_2 = \underline{\Delta}_3 = 0.01$ in the constraints (\ref{eq:NLP:STOConstraints}) to avoid ill-conditioned problems.
The initial guess of the solution was set as $x_0 = x_1 = \cdots = x_N = x (t_0)$, $u_0 = u_1 = \cdots u_{N-1} = 0$, $t_1 = 1.0$, and $t_2 = 2.0$ for all solvers.
It was assumed that the proposed method converged when the solution satisfied the default convergence criteria of Ipopt (max norm of the residual in the unperturbed KKT conditions \cite{bib:nocedal, bib:ipopt}).
All solvers were run 100 times for each $N$ on a laptop with a quad-core CPU Intel Core i7-10510U @1.80 GHz, and the average computational time until convergence was measured.   


\subsubsection{Results}

\begin{table}
\centering
\caption{Average computational time (ms) of each solver for different total number of grids}
\label{table:comparison}
\setlength{\tabcolsep}{3pt}
\begin{tabular}{p{30pt} p{45pt} p{45pt} p{70pt}}
\hline \hline
$N$ & Proposed & {Ipopt} & SQP with {qpOASES} \\ 
\hline
10 & 0.08 & 4.3 & failed \\
50 & 0.27 & 24.6 & failed \\
100 & 0.47 & 45.1 & failed \\
500 & 1.93 & 127 & failed \\
\hline \hline
\end{tabular}
\label{tab1}
\end{table}

The average computational times of each solver for the different total number of grids $N$ are listed in Table \ref{table:comparison}.
As shown in the table, the proposed method converged the fastest of all the cases, and was up to two orders of magnitude faster than Ipopt.
The proposed method was extremely fast because the dimension of the control input of each subsystem was 1 and it did not involve any explicit matrix inversions.
The SQP method with {qpOASES} failed to converge for all cases, even the regularization and indefinite-Hessian mode were enabled.
Note that another QP solver {OSQP} \cite{bib:OSQP} was also tested, but it failed to treat the indefinite Hessian matrix.

The total computational time of Algorithm \ref{alg1} until convergence was also measured, including the mesh-refinement steps.
The thresholds for the discretization step sizes were $\Delta \tau_{\rm max} = 0.35$, $0.065$, $0.035$, and $0.0065$ for cases $N = 10$, $50$, $100$, and $500$, respectively.
The total computational times of Algorithm \ref{alg1} for $N = 10$, $50$, $100$, and $500$, were 0.13, 0.57, 0.92, and 4.3 ms, respectively.
The mesh-refinements were only conducted a few times in each case, and a similar solution was achieved compared to that of the continuous-time counterpart reported in \cite{bib:twoStage}. 
As expected, the solution accuracy improved as $N$ increased.

\subsection{Whole-Body Optimal Control of Quadrupedal Gaits}
\subsubsection{Problem settings}

\begin{figure}[tb]
\centering
\includegraphics[scale=0.29]{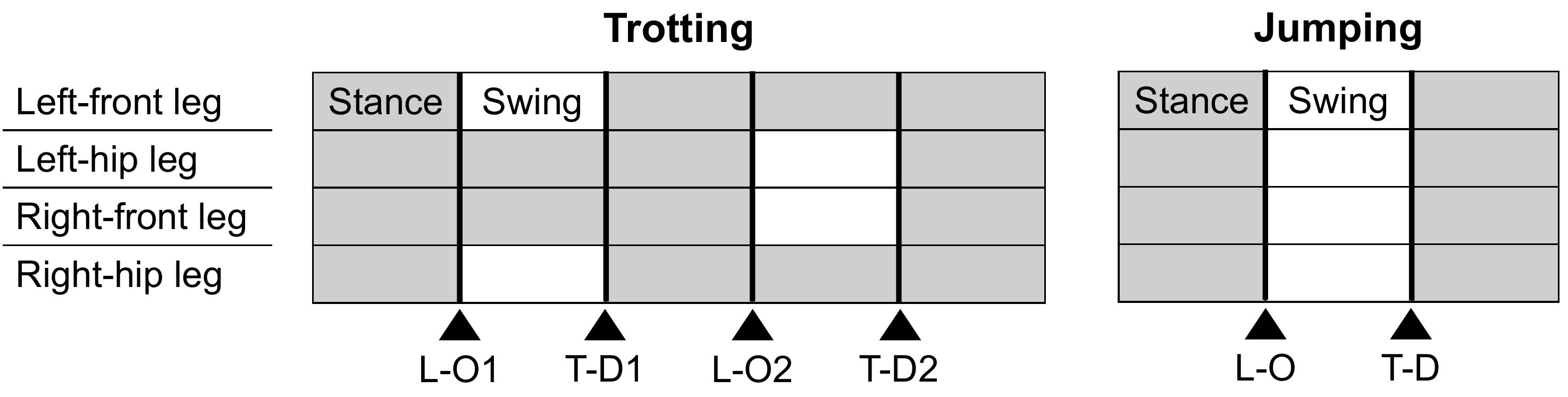}
\caption{Contact patterns of trotting and jumping motions. Gray and white cells indicate the intervals where the leg is standing and swinging, respectively. Thick lines illustrate the switches of the active subsystems. Black triangles indicate the types of switches (lift-off: L-O, touch-down: T-D).}
\label{fig:gaitsequence}
\end{figure}

Numerical experiments on the whole-body optimal control of a quadrupedal robot ANYmal \cite{bib:ANYmal} were conducted to demonstrate the efficiency of the proposed method in practical and complicated examples.
A trotting motion with a step length of 0.15 m and aggressive jumping motion with a jump length of 0.8 m also investigated.
The contact patterns of the trotting and jumping motions are shown in Fig. \ref{fig:gaitsequence}. 
The state equation of the robot switched based on the combination of the support feet, that is, it switched when the feet lifted off from the ground or touched down onto the ground.
Therefore, the trotting motion had four switches and the jump motion had two switches over the horizon, as shown in Fig. \ref{fig:gaitsequence}.
The lift-off did not include the state jump or the switching condition, but the touch-down included both \cite{bib:twoStageApplication1, bib:twoStageApplication2, bib:humanLikeRunnig}.
The switching condition were treated using the proposed method in Section \ref{subseq:switchingConditions}: the switching conditions were position constraints on the foot whose Jacobians took the form of (\ref{eq:pureStateConstraintJacobian}) and the state equation of the robot under an acceleration-level rigid-contact constraint took the form of (\ref{eq:stateEquationForm}). 
The details of the above formulations can be found in \cite{bib:twoStageApplication2, bib:humanLikeRunnig}. 
Note also that the acceleration-level contact-consistent dynamics of the robot were lifted to improve the convergence speed without changing the structure of the KKT system (\ref{subeq:KKTsystem}) \cite{bib:liftedContactDynamics}.
To consider a practical situation, the limitations that were imposed on the joint angles, velocities, and torques were considered to be the inequality constraints.
The polyhedral-approximated friction cone constraint was also considered for each contact force expressed in the world frame $[f_x \;\, f_y \;\, f_z]$ as:
\begin{equation}\label{eq:frictionCone}
    \begin{bmatrix}
        f_x + \frac{\mu}{\sqrt{2}} f_z  \\
        - f_x + \frac{\mu}{\sqrt{2}} f_z  \\
        f_y + \frac{\mu}{\sqrt{2}} f_z  \\
        - f_y + \frac{\mu}{\sqrt{2}} f_z  \\
        f_z
    \end{bmatrix} \geq 0, 
\end{equation}
where $\mu > 0$ is the friction coefficient, which was set as $\mu = 0.7$.
Note that the friction cone constraint (\ref{eq:frictionCone}) also switched depending on the active subsystems (that is, depending on the combination of the support feet). 
The initial and terminal times of the horizon were set as $t_0 = 0$ and $t_{5} = 1.0$ for the trotting problem and $t_0 = 0$ and $t_2 = 1.7$ for the jumping problem, respectively.

To design the stage cost of the trotting motion, the state $x \in \mathbb{R}^{36}$ was partitioned into the generalized coordinate and velocity $[q^{\rm T} \;\; v^{\rm T}] ^{\rm T}$, $q, v \in \mathbb{R}^{18}$. Additionally, the position (i.e., forward kinematics) of the foot $i \in$ \{Left-front, Left-hip, Right-front, Right-hip\} was considered in the world frame $p_i (q) \in \mathbb{R}^3$.
The stage cost was then designed for the trotting motion, which is expressed as:
\begin{align}\label{eq:stageCostTrotting}
    & \frac{1}{2} (q - q_{\rm ref}) ^{\rm T} W_{q} (q - q_{\rm ref})
    + \frac{1}{2} (v - v_{\rm ref}) ^{\rm T} W_{v} (v - v_{\rm ref}) \notag \\ 
    & + \frac{1}{2} a ^{\rm T} W_{a} a \notag \\ 
    & + \sum_{i \in \left\{\rm Swing \; legs \right\}} \frac{1}{2} (p_i (q) - p_{i, {\rm ref}}) ^{\rm T} W_p (p_i (q) - p_{i, {\rm ref}}),
\end{align}
where $a \in \mathbb{R}^{18}$ is the generalized acceleration, $q_{\rm ref}, v_{\rm ref} \in \mathbb{R}^{18}$ and $p_{i, {\rm ref}} \in \mathbb{R}^3$ are reference values, and $W_q, W_v, W_a \in \mathbb{R}^{18 \times 18}$ and $W_p \in \mathbb{R}^{3 \times 3}$ are diagonal weight matrices.
$q_{\rm ref}$ was chosen as the generalized coordinate of the standing pose of the robot.
We set the desired translational velocity of the floating base of the robot as $v_{\rm ref}$ and set the other element of $v_{\rm ref}$ as zero.
The plane coordinate of $p_{i, {\rm ref}}$ was set as a middle point of the two successive predefined ground locations of the foot, and its vertical coordinate was set as the desired height of the swing foot.
The last term in (\ref{eq:stageCostTrotting}) changed depending on the swing legs, that is, the stage cost also switched as the active subsystem.
The terminal cost $V_{f} (x_N)$ and impulse cost $l (x_{i_k -})$ were set as the sums of the first and second terms of (\ref{eq:stageCostTrotting}), respectively.
Similarly, the stage cost was designed for the jumping motions, which is expressed as:
\begin{equation}\label{eq:stageCostJumping}
    \frac{1}{2} (q - q_{\rm ref}) ^{\rm T} W_{q} (q - q_{\rm ref})
    + \frac{1}{2} v ^{\rm T} W_{v} v  
    + \frac{1}{2} a ^{\rm T} W_{a} a . 
\end{equation}
In addition, the terminal cost $V_{f} (x_N)$ and impulse cost $l (x_{i_k -})$ were set as the sums of the first and second terms of (\ref{eq:stageCostJumping}), respectively, but with different weight parameters.
The minimum-dwell times in (\ref{eq:NLP:STOConstraints}) for the trotting motion were set as $\underline{\Delta}_1 = \underline{\Delta}_3 = \underline{\Delta}_5 = 0.02$ and $\underline{\Delta}_2 = \underline{\Delta}_4 = 0.2$.
The minimum-dwell times were set for the jumping motion as $\underline{\Delta}_1 = \underline{\Delta}_2 = 0.15$ and $\underline{\Delta}_3 = 0.65$.
A large $\underline{\Delta}_3$ was set in the jumping motion to sufficiently observe the whole-body motion after touch-down. This was because the optimizer tended to make the touch-down time very close to the end of the horizon to minimize the overall cost of the OCP without a large $\underline{\Delta}_3$. 

The proposed algorithms were implemented in C++ and Eigen was used for the linear algebra. The efficient Pinocchio C++ library was used \cite{bib:pinocchio} for rigid-body dynamics \cite{bib:featherstone} and its analytical derivatives \cite{bib:analyticalRBD} in order to compute the dynamics and its derivatives of the quadrupedal robot.
The Gauss-Newton Hessian approximation was used to avoid computing the second-order derivatives of the rigid-body dynamics.
The reduced Hessian modifications (\ref{eq:Pmodification}) and (\ref{eq:sigmaModification}) were used in the proposed Riccati recursion.
$\sigma_{\rm min}$ and $\bar{\sigma}$ were chosen using (\ref{eq:sigmaModification2}) with $\Delta t_{k, {\rm max}} = 0.1$.
Only the fraction-to-boundary rule \cite{bib:nocedal, bib:ipopt} was used for the step size selection.
The barrier parameter was fixed at $\epsilon = 1.0 \times 10 ^{-3}$, which is a common suboptimal MPC setting \cite{bib:MPCBoyd}.
OpenMP \cite{bib:OpenMP} was used for parallel computing of the KKT system (\ref{subeq:KKTsystem}) in stage-wise and eight threads through the following experiments. 
These two experiments were conducted on a desktop computer with an octa-core CPU Intel Core i9-9900 @3.10 GHz.
$\Delta \tau_{\rm max} = 0.02$ was used for the mesh refinement.
The total number of horizon grids were fixed to $N=50$ and $N=90$ for the trotting and jumping problems, respectively.
That is, when grids were added to the mesh-refinement phase, the same number of grids were removed from the other phases.
The proposed method converged when the $l_2$ norm of the KKT residual was smaller than a sufficiently small threshold $1.0 \times 10 ^{-7}$ and each $\Delta \tau_k$ was smaller than $\Delta \tau_{\rm max}$.
Mesh refinement was performed when the $l_2$ norm of the KKT residual became moderately small (smaller than $0.1$).
This implementation is available online as a part of our efficient optimal control solvers for robotic systems, called robotoc \cite{bib:robotocWeb}.

\subsubsection{Results}
\begin{figure}[tb]
\centering
\includegraphics[scale=0.64]{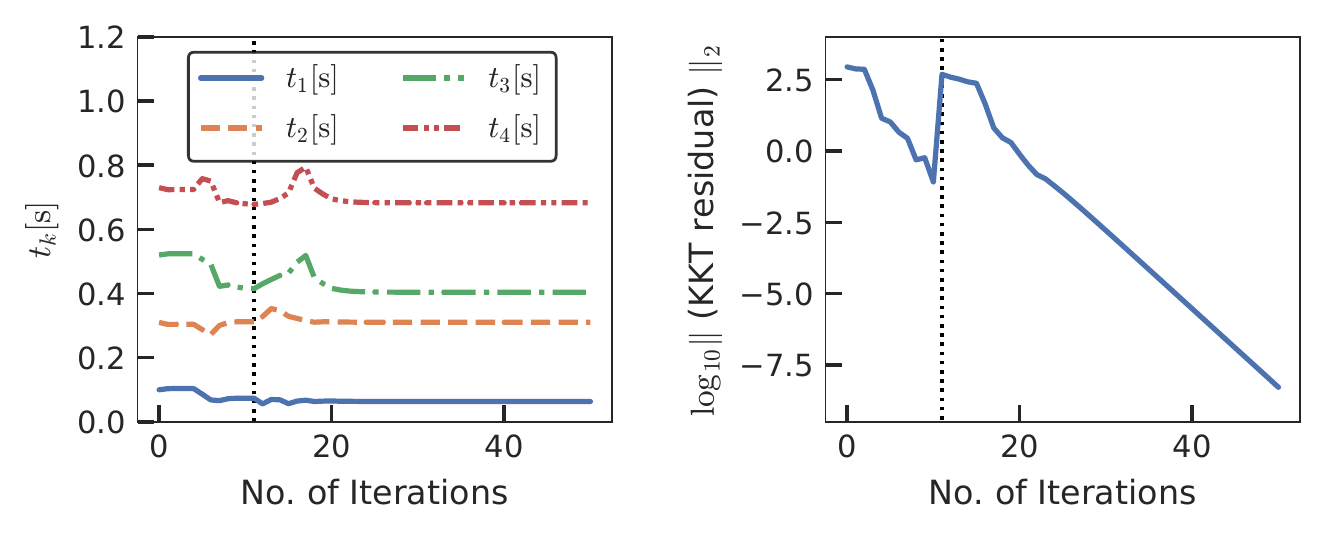}
\caption{Convergence of the proposed method for the trotting motion of a quadrupedal robot: (a) switching instants and (b) squared norm of the residual in the perturbed KKT conditions over the iterations.
Vertical dotted lines indicate that the mesh refinement was carried out.
}
\label{fig:kkt_trot}
\end{figure}

\begin{figure}[tb]
\centering
\includegraphics[scale=0.64]{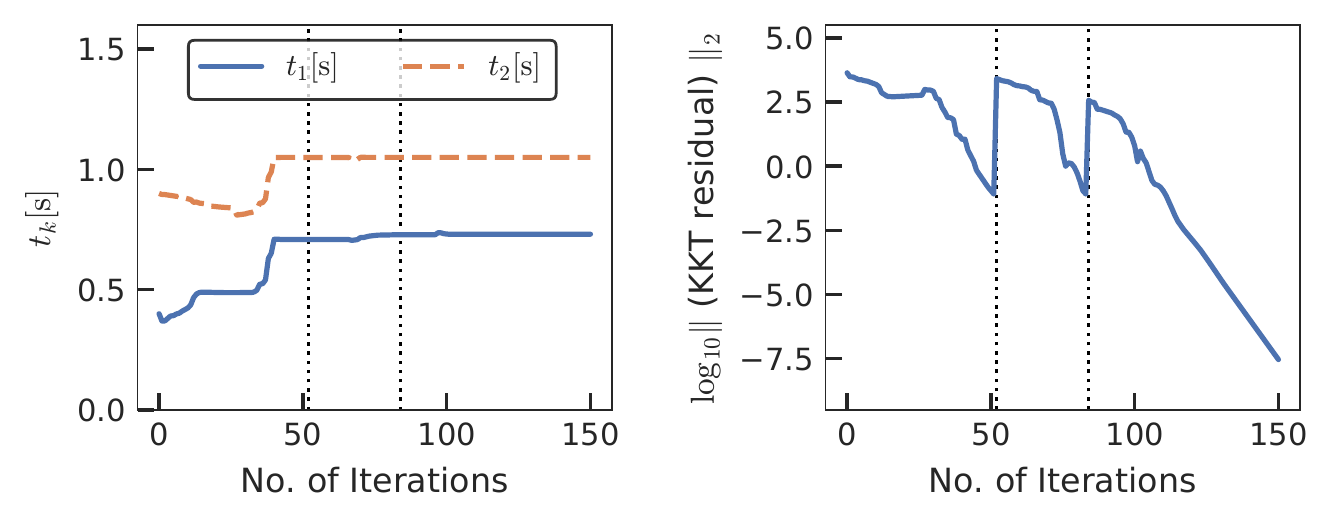}
\caption{Convergence of the proposed method for the jumping motion of the quadrupedal robot: (a) switching instants and (b) squared norm of the residual in the perturbed KKT conditions over the iterations.
Vertical dotted lines indicate that the mesh refinement was carried out.
}
\label{fig:kkt_jump}
\end{figure}

\begin{figure*}[tb]
\centering
\includegraphics[scale=0.06]{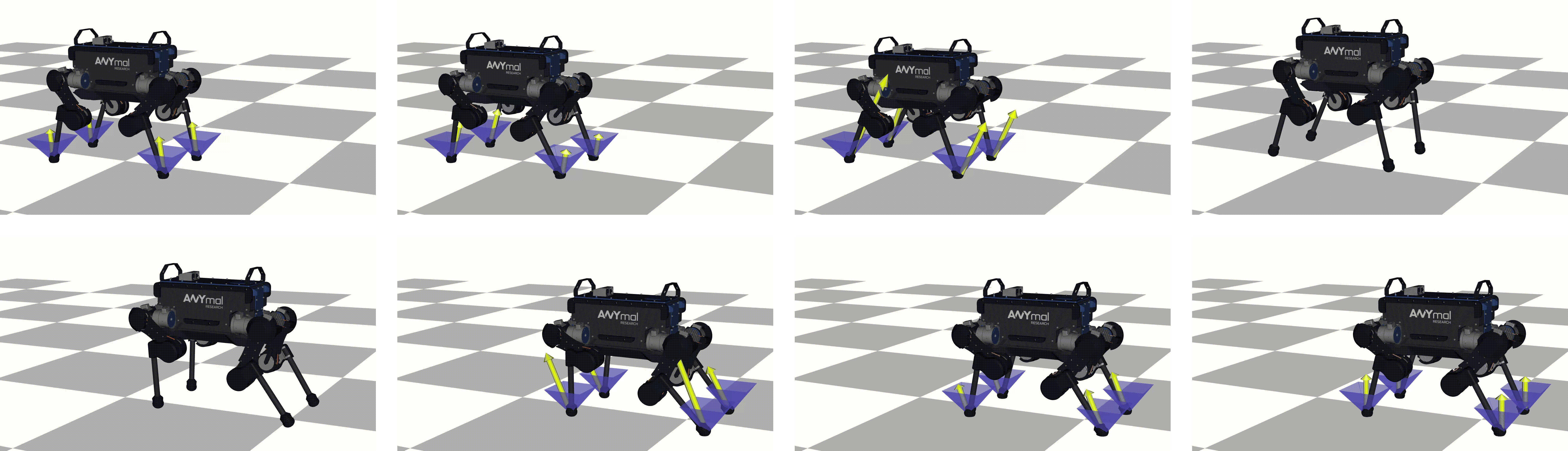}
\caption{Snapshots of solution trajectory of the whole-body optimal control of ANYmal's 0.8 m aggressive jumping motion. The yellow arrows and blue polyhedrons represent the contact forces and friction cone constraints, respectively, which are not illustrated while the robot is flying.
}
\label{fig:jumping}
\end{figure*}

\begin{figure}[tb]
\centering
\includegraphics[scale=0.67]{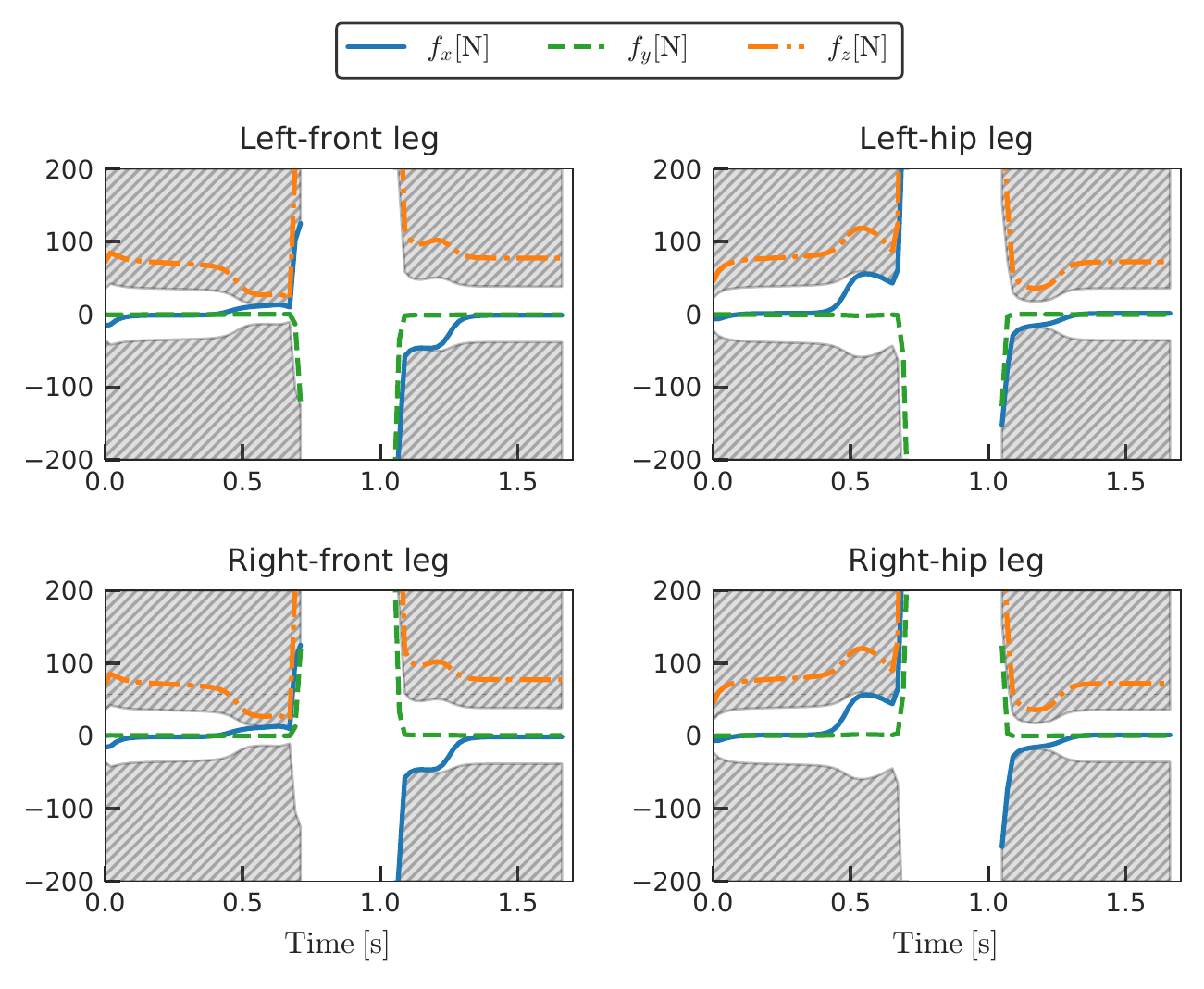}
\caption{Time histories of the contact force expressed in the world frame $[f_x \;\, f_y \;\, f_z]$ of each leg in the jumping motion. The infeasible regions of $f_x$ (solid lines) and $f_y$ (dashed lines) due to the friction cone constraints are the filled gray hatches. 
The infeasible region of $f_z \geq 0$ (dash-dotted lines) is in the lower-half of each plot.}
\label{fig:force}
\end{figure}

The convergence results of the trotting and jumping motions are shown in Figs. \ref{fig:kkt_trot} and \ref{fig:kkt_jump}, respectively.
For the trotting motion that included a mesh-refinement step, the proposed method converged after 46 iterations.
The total computational time was 60 ms (1.3 ms per Newton iteration).
For the jumping motion that included two mesh-refinement steps, the proposed method converged after 147 iterations.
The total computational time was 308 ms (2.1 ms per Newton iteration).
The snapshots of the solution trajectory of the aggressive jumping motion are shown in Fig. \ref{fig:jumping}, and the time histories of the contact forces of the four feet in the jumping motion are shown in Fig. \ref{fig:force}. We can see that the contact forces lie inside the friction cone constraints throughout the jumping motion.
These two figures show that the friction cone constraints (\ref{eq:frictionCone}) were active at time intervals immediately before and after the jumping, and the solution strictly satisfied the constraints.

The above results showed that the proposed method could achieve fast convergence and computational times per Newton iteration, even for large-scale (36-dimensional state, 12-dimensional control input) and complicated problems.
In MPC, the shorter length of the horizon and smaller number of discretization grids can typically be considered, and warm-starting can be leveraged.
For example, it can be expected that the computational time per Newton iteration for the trotting problem would be less than half of 1.3 ms if the number of grids was reduced to 20.
Therefore, the proposed method is a promising approach that can achieve the whole-body MPC of robotics systems with rigid contacts within a milliseconds-range sampling time.

\section{Conclusion}\label{sec:conclu}
This study proposed an efficient Newton-type method for optimal control of switched systems under a given mode sequence.
A direct multiple-shooting method with adaptive mesh refinement was formulated to guarantee the local convergence of the Newton-type method for the NLP.
A dedicated structure-exploiting Riccati recursion algorithm was proposed that performed the Newton-type method for the NLP with the linear time-complexity of the total number of discretization grids.
Moreover, the proposed method could always solve the KKT systems if the solution was sufficiently close to a local minimum. Conversely, general QP solvers cannot be guaranteed to solve the KKT systems because the Hessian matrix is inherently indefinite.
Moreover, to enhance the convergence, a modification on the reduced Hessian matrix was proposed using the nature of the Riccati recursion algorithm as the dynamic programming for a QP subproblem. 
A numerical comparison was conducted with off-the-shelf solvers and demonstrated that the proposed method was up to two orders of magnitude faster.
Whole-body optimal control of quadrupedal gaits was also investigated and it was demonstrated that the proposed method could achieve the whole-body MPC of robotic systems with rigid contacts.

A possible future extension of the proposed method is OCPs with free-final time, including minimum-time OCPs \cite{bib:appliedOCP}, because the NLP structure of such problems is expected to be similar to the proposed formulation.

\section*{Acknowledgment}
This work was partly supported by JST SPRING, Grant Number JPMJSP2110.

\bibliographystyle{ieeecolor}
\bibliography{IEEEabrv, root}

\end{document}